\documentclass[11pt,letterpaper]{article}
\usepackage[in]{fullpage}

\usepackage[citestyle=authoryear,bibstyle=authoryear, natbib=true, backend=bibtex,maxnames=2, minnames=1, maxbibnames=99, giveninits, uniquename=false]{biblatex}
\newrobustcmd*{\citefirstlastauthor}{\AtNextCite{\DeclareNameAlias{labelname}{given-family}}\citeauthor*}
\DeclareFieldFormat*{title}{``#1''}

\usepackage[utf8]{inputenc}
\usepackage{amsmath}
\usepackage{amsfonts}
\usepackage{amssymb}
\usepackage{amsthm}
\usepackage{dsfont}
\usepackage{xcolor}
\usepackage{hyperref}
\usepackage{mathtools}
\usepackage{pgf}
\usepackage{chngpage}
\usepackage{subfigure}
\usepackage{longtable}
\usepackage[noabbrev,capitalize]{cleveref}

\DeclarePairedDelimiter{\braces}{\{}{\}}
\DeclarePairedDelimiter{\bracks}{[}{]}
\DeclarePairedDelimiter{\parens}{(}{)}

\DeclarePairedDelimiterX{\braket}[2]{\langle}{\rangle}{#1,#2}

\DeclarePairedDelimiterX{\inner}[2]{\langle}{\rangle}{#1,#2}
\DeclarePairedDelimiterX{\setdef}[2]{\{}{\}}{#1:#2}

\DeclarePairedDelimiterXPP{\probof}[1]{\prob}{(}{)}{}{%

#1}

\DeclarePairedDelimiterXPP{\exof}[1]{\ex}{[}{]}{}{%

#1}

\newcommand{\argdot}{\,\cdot\,}

\DeclareMathOperator*{\argmin}{\arg\min}

\usepackage{bbding}
\usepackage{pifont}
\newcommand{\cmark}{\ding{51}}%

\usepackage{array}
\newcolumntype{x}[1]{>{\centering\arraybackslash\hspace{0pt}}p{#1}}

\newtheorem{assumption}{Assumption}
\newtheorem{theorem}{Theorem}
\newtheorem{lemma}{Lemma}
\numberwithin{equation}{section}
\numberwithin{theorem}{section}
\numberwithin{lemma}{section}

\crefname{assumption}{Assumption}{Assumptions}

\newcommand{\ie}{i.e.,\ }
\newcommand{\eg}{e.g.,\ }

\newcommand{\one}{\mathds{1}}
\newcommand{\sign}[1]{\mathrm{sign}\left(#1\right)}
\newcommand{\iid}{i.i.d.\ }

\newcommand{\loss}{\mathcal{L}_T}
\newcommand{\regret}{R_T}
\newcommand{\diffbound}{M}
\newcommand{\imp}{\mathrm{naive}}

\let\oldmaketitle\maketitle
\renewcommand{\maketitle}{\oldmaketitle\setcounter{footnote}{0}}
\newcommand{\munderbar}[1]{\text{\underline{$#1$}}}

\newcommand{\bigO}[1]{\mathcal{O}\left(#1\right)}
\newcommand{\kl}[1]{\mathrm{KL}\left(#1\right)}

\makeatletter
\newcommand*\rel@kern[1]{\kern#1\dimexpr\macc@kerna}
\newcommand*\widebar[1]{%
  \begingroup
  \def\mathaccent##1##2{%
    \rel@kern{0.8}%
    \overline{\rel@kern{-0.8}\macc@nucleus\rel@kern{0.2}}%
    \rel@kern{-0.2}%
  }%
  \macc@depth\@ne
  \let\math@bgroup\@empty \let\math@egroup\macc@set@skewchar
  \mathsurround\z@ \frozen@everymath{\mathgroup\macc@group\relax}%
  \macc@set@skewchar\relax
  \let\mathaccentV\macc@nested@a
  \macc@nested@a\relax111{#1}%
  \endgroup
}

\title{Social Learning in non-stationary environments}

\author{Etienne Boursier \\
TML Lab, EPFL, Lausanne, Switzerland\\
etienne.boursier1@gmail.fr
 \and Vianney Perchet \\
 CREST, ENSAE Paris, Palaiseau, France  \\
 CRITEO AI Lab, Paris, France\\
 vianney.perchet@normalesup.or
  \and Marco Scarsini\\
  LUISS University, Rome, Italy\\
 marco.scarsini@luiss.it
  }
\date{}

\begin{document}
\maketitle              

\begin{abstract}

Potential buyers of a product or service, before making their decisions, tend to read reviews written by previous consumers. 
We consider Bayesian consumers with heterogeneous preferences, who sequentially decide whether to buy an item of unknown quality, based on previous buyers' reviews. 
The quality is multi-dimensional and may occasionally vary over time; the reviews are also multi-dimensional.
In the simple uni-dimensional and static setting, beliefs about the quality are known to converge to its true value. 
Our paper extends this result in several ways. First, a multi-dimensional quality is considered, second, rates of convergence are provided, third, a dynamical Markovian model with varying quality is studied.  
In this dynamical setting the cost of learning is shown to be small. 


\end{abstract}

\section{Introduction}
\label{se:intro}

In our society many forms of learning do not stem from direct experience, but rather from observing the behavior of other people who themselves are trying to learn. 
In other words, people engage in social learning. 
For instance, before deciding whether to buy a product or service, consumers observe the past behavior of previous consumers and use this observation to make their own decision. 
Once their decision is made, this becomes a piece of information for future consumers.
In the old days, it was common to consider a crowd in a restaurant as a sign that the food was likely good. 
Nowadays, there are more sophisticated ways to learn from previous consumers. 
After buying a product and experiencing its features, people often leave reviews on sites such as Amazon, Tripadvisor, Yelp, etc.
When consumers observe only the purchasing behavior of previous consumers, there is a risk of a cascade of bad decisions: if the first agents make the wrong decision, the subsequent agents may follow them thinking that what they did was optimal and herding happens. 
Interestingly enough, this is not necessarily the effect of bounded rationality.
It can actually be the outcome of a Bayesian equilibrium in a game with fully rational players.
It seems reasonable to conjecture that, if consumers write reviews about the product that they bought, then social learning will be achieved. 
This is not always the case when consumers are heterogeneous and the reviews that they write depend on the quality of the object but also on their idiosyncratic attitude towards the product they bought. 

Consumers also tend to give higher value to recent reviews. 
As highlighted in a survey \citep{murphy2019local} run on a panel of a thousand consumers,  ``48\% of consumers only pay attention to reviews written within the past two weeks,'' and this percentage is increasing over the years.
A justification for this behavior may be that customers perceive the quality of the product that they consider buying as variable over time. 
The more recent the review, the more informative it is about the current state of the product. 
This work considers a dynamical environment and shows that, under some conditions,  the outcome of the learning process in  stationary and non-stationary environments are overall comparable.

\subsection{Main contribution}
\label{suse:contribution}

We consider a model where heterogeneous consumers arrive sequentially at a monopolistic market and---before deciding whether to buy a product of unknown quality---observe the reviews (\eg like/dislike) provided by previous buyers. 
Consumers are Bayesian and buy the product if and only if their expected utility of buying is larger than $0$ (the utility of the outside option).
Each buyer posts a sincere review that summarizes the experienced quality of the product, affected by an idiosyncratic attitude to the product. 
\citet{ifrach2019} studied this model in the case where the intrinsic quality of the product is one-dimensional, fixed over time, and can assume just two values; they studied conditions for social learning to be achieved. 
We extend their results in two main directions. 
First, we allow the quality to be multidimensional, \ie to have different features that consumers experience and evaluate.
Second, we consider a model where the quality can occasionally change over time.

We start examining a benchmark model where the quality is actually static and we provide  rates of convergence for the posterior distribution of the quality.
We then move to the more challenging dynamical model where quality may change over time. 
The criterion that we use in this dynamical setting is the utility loss that a non-informed consumer incurs with respect to a fully informed consumer, who at every time knows the true quality of the product.
We show that  the learning cost is a logarithmic factor of the changing rate of the quality.

\cref{table:summary} below summarizes the proved bounds for the different settings. 
In the analysis we also consider the case of imperfect learners, who are not aware of the dynamical nature of the quality, and we quantify the loss they incur.

\begin{table}[ht]
\centering\setlength{\extrarowheight}{4pt}
\setlength\tabcolsep{5pt}

\begin{tabular}{|c||c|c|}
\hline
\textbf{Type of model}  & \textbf{Utility Loss} & \textbf{Tight Bound}\\[4pt]
\hline
\hline
stationary &   $\bigO{Md}$ & \cmark \\[4pt]
\hline
dynamical &   $\bigO{Md\ln(2/\eta)\eta T}$ & \cmark \\[4pt]
\hline
\end{tabular}
\vspace{0.5em}
\caption{\label{table:summary}Bounds summary, where the reward function is $M$-Lipschitz and $d$ is the dimension of the quality space. In a non-stationary environment, the quality changes with probability $\eta$ at each round, while the utility loss is summed over $T$ rounds.}
\vspace{-2em}
\end{table}

%
%

\subsection{Related literature}
\label{suse:related-literature}

The problem of social learning goes back to \citet{banerjee1992simple,bikhchandani1992theory} who considered models where Bayesian rational agents arrive at a market sequentially, observe the actions of the previous agents, and decide based on their private signals and the public observations.
These authors showed that in equilibrium, consumers may herd into a sequence of bad decisions; in other words, social learning  fails with positive probability.
\citet{smith2000pathological} showed that this learning failure is due to the fact that signals are bounded. 
In the presence of unbounded signals that can overcome any observed behavior, herding cannot happen.

Different variations of the above model have been considered, where either agents observe only a subset of the previous agents \citep[see \eg][]{ccelen2004observational,acemoglu2011bayesian,lobel2015information}, or the order in which actions are taken is not determined by a line, but rather by a lattice \citep{arieli2019multi}.
A general analysis of social learning models can be found in \citep{AriMue:MOR2021}.

A more recent stream of literature deals with models where agents observe not just the actions of the previous agents, but also their ex-post reaction to the actions they took. 
For instance, before buying a product of unknown quality, consumers read the reviews written by the previous consumers. 
In particular, \citet{besbes2018} dealt with some variation of a model of social learning in the presence of reviews with heterogeneous consumers. 
In one case, agents observe the whole history of reviews and can use Bayes rule to compute the conditional expectation of the unknown quality and learning is achieved.
In the other case they only observe the mean of past reviews.
Interestingly, even in this case, learning is achieved and the speed of convergence is of the same order.
\citet{ifrach2019} studied a model where the unknown quality is binary and the reviews are also binary (like or dislike). 
They considered the optimal pricing policy and looked at conditions that guarantee social learning. 
\citet{correa2020} also considered the optimal dynamic pricing policy when consumers have homogeneous preferences.
A non-Bayesian version of the model was considered by \citet{crapis2016monopoly}, where mean-field techniques were adopted to study the learning trajectory.

\citet{papanastasiou2017dynamic} studied a market where strategic consumers can delay their purchase anticipating the fact that other consumers will write reviews in the meanwhile. 
They examined the implication on pricing of this strategic interaction between consumers and a monopolist. 
\citet{feldman2018social} examined the role of social learning from reviews in the monopolist's design of a product in a market with strategic consumers.
\citet{KakLan:MIT2021} studied heterogeneity of consumers' reviews and its impact on social learning and price competition.
\citet{MagScaVac:SSRN2020} considered a model of social learning with reviews where consumers have different buying options and a platform can affect consumers' choice by deciding the order in which different brands are displayed.
\citet{ParShiXie:MS2021} dealt with the effect of the first review on the long-lasting success of a product.
\citet{CheLiTal:MS2021} considered the issue of bias in reviews from a theoretical viewpoint. They quantified the acquisition bias and the impact on the rating of an arriving customer, characterized the asymptotic outcome of social learning, and we show the effect of biases and social learning on pricing decisions. 

The speed of convergence in social learning was considered by \citet{rosenberg2019} in  models where only the actions of the previous agents are observed and by \citet{acemoglu2017} when reviews are present. 
This last paper is the closest to the spirit of our own paper.

Learning problems in non-stationary environment have been largely considered in online learning with bandit feedback. Two different approaches are generally considered. First, the environment can be abruptly changing, either following a non-parametric model \citep{yu2009piecewise,garivier2011upper,besson2019efficient} or a Markov chain \citep{mellor2013thompson,alami2018memory} as in our model in \cref{sec:dynamic}. The regret then scales with the number of changes in the environment parameters.
Conversely, the environment can be changing smoothly with time \citep[see \eg][]{besbes2015nonstationary,besbes2019optimal,keskin2017chasing,cheung2018hedging}. The regret then scales with the cumulated change in the parameters $B = \sum_{t=1}^{T-1} \|\theta_{t+1} - \theta_t\|$, where $\theta_t$ is the parameter vector at time~$t$. These models thus differ both in their practical motivations and technical solutions. We consider in this work the former setting, which is natural to use within a Bayesian framework.

Social learning however differs from bandits models, as the users act greedily and only maximize their reward over a single round. In bandits, a single meta-user chooses over repeated steps and thus maximizes its reward over multiple steps. 
Works in social learning often depict conditions under which the greedy policy still leads to learning the parameters, while it is unknown to generally fail in bandits problems.

\subsection{Organization of the paper}
\label{suse:organization}

\cref{sec:model} introduces the model of social learning from consumer reviews. 
\cref{sec:stationary} studies the stationary setting where the quality is fixed.
\cref{sec:dynamic} introduces the dynamical setting, where the quality changes over time. 
\cref{sec:imperfect_learners} consider a model with naive consumers and shows that knowledge of the dynamical structure is crucial for the consumer utility.

\cref{se:symbols} provides a list of symbols, \cref{se:additional-proofs} contains additional proofs. \cref{app:continuous} studies the continuous model where the quality space $\mathcal{Q}$ is convex.

%
%

\section{Model}\label{sec:model}


We consider a model of social learning where consumers read reviews before making their purchase decisions. 
A monopolist sells a product of unknown quality to consumers who arrive sequentially at the market. 
The quality may vary over time, although variations are typically rare.
The quality of the product at time $t$ is denoted by $Q_{t}$ and
the set of possible qualities is $\mathcal{Q}=\lbrace 0, 1 \rbrace^d$. For a vector $x$, we denote by $x^{(i)}$ its $i$-th component, \ie $Q_t^{(i)}$ represents the $i$-th feature of the product at time $t$ and has a binary value (low or high).

The prior distribution of the quality at time $1$ is $\pi_{1}$.
Consumers are indexed by the time of their arrivals $t\in\mathbb{N}\setminus\lbrace0\rbrace$.
They are heterogeneous and consumer $t$ has an idiosyncratic preference $\theta_{t} \in \Theta$ for the product. 
This preference $\theta_{t}$ is private information. 
These preferences are assumed to be \iid according to some known distribution.
In game-theoretic terms, $\theta_{t}$ could be seen as the \emph{type} of consumer~$t$.
The sequences of preferences $\theta_{t}$ and of qualities $Q_{t}$ are independent.

A consumer who buys the product posts a review in the form of a  multi-dimensional numerical grade.
The symbol $Z_{t}$ denotes the review posted by consumer $t$.
The notation $Z_{t}=\ast$ indicates that consumer $t$ did not buy the product.
We call $\mathcal{H}_{t} \coloneqq \braces*{Z_{1},\dots,Z_{t-1}}$ the history before the decision of consumer~$t$.
We set $\mathcal{H}_{1} \coloneqq \varnothing$.

Since the preferences are independent of the quality, a no-purchase decision does not carry any information on the quality.
As a consequence, the history $\mathcal{H}_{t}$ is informationally equivalent to the reduced history $\widetilde{\mathcal{H}}_{t}$ that includes only the reviews of the buyers up to $t-1$.
This differentiates this model from the classical social learning models, where consumers have private signals that are correlated with the quality.

%
%
%
%
%


Based on the history $\mathcal{H}_t$ of past observations and her own preference $\theta_t$, consumer~$t$ decides whether to buy the product. 
In case of purchase, she receives the utility $u_{t} \coloneqq r(Q_t, \theta_t)$ where $r$ is the reward function.
A consumer who does not buy the product gets $u_{t} = 0$. 

Bayesian rationality is assumed, so consumer $t$ buys the product if and only if her conditional expected utility of purchasing is  positive, that is, if and only if $\mathbb{E}[r(Q_t, \theta_t) \mid \mathcal{H}_t, \theta_t] > 0$. 
Consumer $t$ then reviews the product by giving the feedback $Z_t = f(Q_t, \theta_t, \varepsilon_t) \in \mathcal{Z} \subset \mathbb{R}^d$ where $\varepsilon_t$ are i.i.d. variables independent from $\theta_t$. 
Also, the feedback function is assumed to take a finite number of values in $\mathbb{R}^d$ and to be of the form
\begin{equation*}
f(Q, \theta, \varepsilon) = (f^{(i)}(Q^{(i)}, \varepsilon, \theta))_{i=1,\ldots,d}\ .
\end{equation*}
In words, for each different feature $Q^{(i)}$ of the quality $Q$, consumers provide a separate feedback. Previous works \citep{acemoglu2017,ifrach2019} considered $\mathcal{Z} = \lbrace 0, 1 \rbrace$  as the reviews were only the likes or dislikes of consumers. This model allows a more general and richer feedback, such as ratings on a five-star scale for each feature, or even sparse feedback where consumers do not necessarily review each feature.
%

In a model without noise $\varepsilon_t$, the learning process is much simpler, as already noted by \citet{ifrach2019}. 
Indeed, in this case, a single negative review rules out many possibilities as it means that the quality was overestimated. 
We here consider noise, which corresponds to variations caused by different factors such as fluctuations in the product quality or imperfect perception of the quality by the consumer. This also leads to a more interesting learning process with respect to the noiseless case.

\medskip

In the following, $\pi_t$ denotes the posterior distribution of $Q_t$ given $\mathcal{H}_t$ and, for any $i\in[d] \coloneqq \lbrace 1, \ldots, d \rbrace$, $\pi_t^{(i)}(q^{(i)}) = \mathbb{P}[Q_t^{(i)} = q^{(i)} \mid \mathcal{H}_t]$ is the $i$-th marginal of the posterior.
%
%
We also introduce the function $G$ and its componentwise equivalent $G^{(i)}$, defined as 
\begin{align}
\label{eq:defnG} 
G(z, \pi, q) &= \mathbb{P}[Z_t=z \mid \pi_t=\pi, Q_t=q],\\
\label{eq:defnGi} 
G^{(i)}(z^{(i)}, \pi, q^{(i)}) &= \mathbb{P}[Z_{t}^{(i)}=z^{(i)} \mid \pi_t=\pi, Q_t^{(i)}=q^{(i)}].
\end{align}
%
In the following, we also use the notations 
\begin{align}
\label{eq:defnGzpi}
G(z, \pi) &= \mathbb{E}_{q \sim \pi}[G(z, \pi, q)],\\
\label{eq:defnGizpi}
G^{(i)}(z^{(i)}, \pi) &= \mathbb{E}_{q \sim \pi}[G^{(i)}(z^{(i)}, \pi, q^{(i)})].
\end{align} 
%
%
The following two assumptions will be used in the sequel.
\begin{assumption}[Purchase guarantee]\label{ass:buyers}
The reward function $r$ is weakly increasing in each feature $q^{(i)}$ and for any $q\in \mathcal{Q}$, $\mathbb{P}_{\theta_t}\big(r(q, \theta_t) > 0 \big) > 0$, \ie there is always a fraction of consumers who buy the product.
\end{assumption}
\cref{ass:buyers} excludes situations where all consumers stop buying if the quality belief becomes too low. 
Without this condition, social learning fails with positive probability \citep{acemoglu2017,ifrach2019}. 
%
\begin{assumption}[Identifiability]\label{ass:identification}
For any $i \in [d]$, any quality posterior $\pi \in \mathcal{P}(\mathcal{Q})$ and quality $q^{(i)}$, we have $G^{(i)}(\argdot, \pi, q^{(i)}) > 0$.
Moreover, for $q^{(i)} \neq q'^{(i)}$, there exists some $z \in \mathcal{Z}$ such that $G^{(i)}(z^{(i)}, \pi, q^{(i)}) \neq G^{(i)}(z^{(i)}, \pi, q'^{(i)})$.
\end{assumption}
\cref{ass:identification} is needed to distinguish different qualities based on past reviews. 
Non-identifiability could lead to a considerable utility loss, as the set of consumers who should buy might differ for two undistinguishable qualities.
The positivity of $G$ is required to avoid trivial situations.  
The case of $G=0$ for some variables is similar to the absence of noise $\varepsilon_t$, as a single observation can definitely rule out several possibilities.
%

%
\medskip

An interesting choice of reward function is, for instance, $r(Q, \theta) = \langle Q, \theta_{t} \rangle$ where $\langle \argdot, \argdot \rangle$ is the scalar product. 
In this case, $\theta_{t}^{(i)}$ is the weight that customer~$t$ gives to  feature $i$ of the service.
In practice, customers might also only focus on the best or worst aspects of the service, meaning their reward might only depend on the maximal or minimal value of the $Q^{(i)}$'s.
The ordered weighted averaging operators \citep{yager1988}  model these behaviors. 
In an additive model similar to the classical case in the literature, this leads to a reward function $r(Q, \theta) = \sum_{i=1}^n w^{(i)} (Q+\theta)^{(\sigma(i))}$ where $\sigma$ is a permutation such that $(Q+\theta)^{(\sigma(i))}$ is the $i$-th largest component of the vector $(Q^{(i)} + \theta^{(i)})_{i=1, \ldots, d}$. 
If $w^{(i)}=1/d$ for all $i$, this is just an average of  all features' utilities. 
When $w^{(1)}=1$ and all other terms are $0$, consumers are only interested in the maximal utility among all  features.

Much of the existing literature has focused on the following unidimensional setting
\begin{align*}
r(Q, \theta) &= Q + \theta - p, \\
f(Q, \theta, \varepsilon) &= \sign{Q + \theta + \varepsilon - p},
\end{align*}
where $p$ is an exogenously fixed price. 
%
%
Since consumers review separately each feature of the service, the feedback function is a direct extension of the above unidimensional setting. 
It is then of the form 
\begin{equation}
\label{eq:fi-sign}
f^{(i)}(Q^{(i)}, \varepsilon, \theta) = \sign{Q^{(i)} + \theta^{(i)} + \varepsilon^{(i)} - p^{(i)}},
\end{equation}
for some constant price $p^{(i)}$. 

Having a sparse feedback is very common on platform reviews, where consumers only review a few features. 
This case  can be modeled by 
\begin{equation}
\label{eq:fi-sparse}
f^{(i)}(Q^{(i)}, \varepsilon, \eta, \theta) = \sign{Q^{(i)} + \theta^{(i)} + \varepsilon^{(i)} - p^{(i)}} \, \xi^{(i)},
\end{equation}
with $\varepsilon \in \mathbb{R}^{d}$ and $\xi \in \lbrace0, 1 \rbrace^d$. 
Although the noise vector is here given by the tuple $(\varepsilon, \xi)$ instead of $\varepsilon$ alone, this remains a specific case of our model.

A multiplicative model can also be considered where the relevant quantity is $Q^{(i)} \theta^{(i)}$, rather than $Q^{(i)} + \theta^{(i)}$. 
This model is very similar to the additive one when using a logarithmic transformation.

\section{Stationary Environment}
\label{sec:stationary}

As mentioned before, our aim is to consider a model where the quality of the product may occasionally change over time. 
As a benchmark, we start considering the case where the quality is constant: $Q_t = Q_1$ for all $t \in \mathbb{N}$. 
We will leverage this case, when dealing with the dynamic model of variable quality.
In the unidimensional case $\mathcal{Q}=\lbrace 0, 1 \rbrace$, \citet{ifrach2019} showed that the posterior almost surely converges to the true quality, and \citet{acemoglu2017} showed an asymptotic exponential convergence rate. 
Besides extending these results to the multidimensional model, this section shows anytime  convergence rates of the posterior. 
Convergence rates in social learning has been studied only recently \citep{acemoglu2017,rosenberg2019}, despite being central to online learning \citep{bottou1999} and Bayesian estimation \citep{ghosal2000}. 
Moreover, convergence rates are of crucial interest when facing a dynamical quality. The main goal of this section is thus to lay the foundation for the analysis of \cref{sec:dynamic}.

The posterior update is obtained using Bayes' rule for any $q\in \mathcal{Q}$,
\begin{equation}\label{eq:updatestationarydiscretefull}
\pi_{t+1}(q) = \frac{G\left(Z_t,\pi_t, q\right)}{G\left(Z_t,\pi_t \right)} \ \pi_t(q).
\end{equation}
\cref{thm:stationary_discrete} below gives a convergence rate of the posterior to the true quality. Similarly to  \citet[Theorem~2]{acemoglu2017}, it shows an exponential convergence rate. 
The main difference with \citet[Theorem~2]{acemoglu2017} is that we here provide an anytime bound, while their convergence rate is only asymptotic, \ie they show that almost surely $\lim_{t\to \infty} log(pi_t(q’))/t = -c$ for some constant $c$ of the same order than the exponential rate in \cref{thm:stationary_discrete}.
We focus on anytime rates as they are highly relevant in the model with a dynamical, evolving quality considered in \cref{sec:dynamic}.
\begin{theorem}\label{thm:stationary_discrete}
For $q \neq q'$, we have 
\begin{equation*}
\mathbb{E}[\pi_{t+1}(q') \mid Q=q] \leq \exp\left(-\frac{t \delta^4}{2\gamma^2 + 4 \delta^2}\right) \frac{1}{\max_{i\in[d]} \pi_{1}^{(i)}(q^{(i)})},
\end{equation*}
\begin{align}
&\text{where} \quad \label{eq:delta-i}
\delta \coloneqq \min_{\substack{i\in [d], \pi \in \mathcal{P}(\mathcal{Q})}} \sum_{z \in \mathcal{Z}}| G^{(i)}(z^{(i)}, \pi, 1) - G^{(i)}(z^{(i)}, \pi, 0)| \\
\label{eq:gamma-i}
&\text{and} \quad \gamma \coloneqq 2\max_{\substack{i\in [d], \pi \in \mathcal{P}(\mathcal{Q}), z \in \mathcal{Z}}}  \left| \ln\left(\frac{G^{(i)}(z^{(i)}, \pi, 1)}{G^{(i)}(z^{(i)}, \pi, 0)}\right) \right|.
\end{align} 
\end{theorem}
Notice that $\delta$ is the minimal total variation between $Z_t^{(i)}$ conditioned either on ${(\pi, Q_t^{(i)}=1)}$ or ${(\pi, Q_t^{(i)}=0)}$. Thanks to \cref{ass:identification}, both $\delta$ and $\gamma$ are positive and finite. 
This guarantees an exponential convergence rate of the posterior as $\pi_t(q) = 1 - \sum_{q' \neq q} \pi_t(q')$.

\begin{proof}
Assume without loss of generality $Q_1^{(i)}=1$. 
\cref{thm:stationary_discrete} is a direct consequence of \cref{eq:marginalconvergence}, which we prove in the following:
\begin{equation}\label{eq:marginalconvergence}
\mathbb{E}[\pi_{t+1}^{(i)}(0) \mid Q_1^{(i)}=1] \leq \exp\left(-\frac{t \delta^4}{2\gamma^2 + 4 \delta^2}\right) \frac{1}{\pi_1^{(i)}(1)}.
\end{equation}

Similarly to \cref{eq:updatestationarydiscretefull}, we have the Bayesian update
$
\pi_{t+1}^{(i)}(q^{(i)}) = \frac{G^{(i)}\left(Z_t^{(i)},\pi_t, q^{(i)}\right)}{G^{(i)}\left(Z_t^{(i)},\pi_t \right)} \ \pi_t^{(i)}(q^{(i)})$, which leads by induction to
\begin{equation*}
\ln\left(\frac{\pi^{(i)}_{t+1}(1)}{\pi^{(i)}_{t+1}(0)}\right) =\ln\left(\frac{\pi^{(i)}_{1}(1)}{\pi^{(i)}_{1}(0)}\right) + \sum_{s=1}^t \ln\left(\frac{G^{(i)}\left(Z_s^{(i)},\pi_s, 1\right)}{G^{(i)}\left(Z_s^{(i)}, \pi_s, 0\right)}\right).
\end{equation*} 

In the following, we  use the notation $\kl{\mu, \nu}$ for the Kullback-Leibler divergence between the distributions $\mu$ and $\nu$, which is defined as
\begin{equation}
\label{eq:KL}
\kl{\mu,\nu} = \mathbb{E}_{x \sim \mu} \left[ \ln\left(\frac{\mu(x)}{\nu(x)}\right)\right].
\end{equation}
Now define $
X_t \coloneqq \ln\left(\frac{G^{(i)}(Z_t^{(i)}, \pi_t, 1)}{G^{(i)}(Z_t^{(i)}, \pi_t, 0)} \right) - \kl{G^{(i)}(\cdot, \pi_t, 1), G^{(i)}(\cdot, \pi_t, 0)}$.

%
Notice that $\mathbb{E}[X_t \mid \mathcal{H}_t, Q_1^{(i)}=1] = 0$. Also, by definition of $\gamma$, $X_t \in [Y_t, Y_t+\gamma]$ almost surely for some $\mathcal{H}_t$-measurable variable $Y_t$.
Azuma-Hoeffding's inequality for martingales \citep[see, e.g.,][Lemma~A.7]{cesa2006} then yields for any $\lambda \geq 0$:
\begin{align*}
\mathbb{P}\left[\sum_{s=1}^t X_s \leq -\lambda \,\Big|\, Q_1^{(i)}=1\right] \leq \exp\left(-\frac{2\lambda^2}{t\gamma^2}\right),
\end{align*}
which is equivalent to
\begin{equation*}
\begin{split}
&\mathbb{P}\left[\frac{\pi^{(i)}_{t+1}(0)}{\pi^{(i)}_{t+1}(1)} \geq \exp\Big(\lambda - \sum_{s=1}^t \kl{G^{(i)}(\cdot, \pi_s, 1), G^{(i)}(\cdot, \pi_s, 0)}\Big)\frac{\pi_1^{(i)}(0)}{\pi_1^{(i)}(1)}  \,\Big|\, Q^{(i)}_1=1 \right] \leq \exp\left(-\frac{2\lambda^2}{t\gamma^2}\right).
\end{split}
\end{equation*}
By Pinsker's inequality \citep[see, \eg][Lemma~2.5]{Tsy:Springer2009}, we have
\begin{equation*}
\kl{G^{(i)}(\cdot, \pi_s, 1), G^{(i)}(\cdot, \pi_s, 0)} \geq \delta^2/2,
\end{equation*}
so the previous equation becomes
\begin{equation*}
\mathbb{P}\left[\pi^{(i)}_{t+1}(0) \geq \exp\Big(\lambda - \frac{t\delta^2(1, 0)}{2} \Big)\frac{\pi_1^{(i)}(0)}{\pi_1^{(i)}(1)} \,\Big|\, Q_1^{(i)}=1 \right] \leq \exp\left(-\frac{2\lambda^2}{t\gamma^2}\right),
\end{equation*}
where we used the fact that $\pi^{(i)}_{t+1}(1) \leq 1$. 
This then yields
\begin{align*}
\mathbb{E}[\pi^{(i)}_{t+1}(0) \mid Q_1^{(i)}=1] & \leq \exp\left(\lambda-\frac{t\delta^2}{2}\right)\frac{\pi_1^{(i)}(0)}{\pi_1^{(i)}(1)} + \mathbb{P}\left[\pi^{(i)}_{t+1}(0) \geq \exp\left(\lambda - t\delta^2/2 \right) \mid Q_1^{(i)}=1 \right]\\
& \leq \exp\left(\lambda-\frac{t\delta^2}{2}\right)\frac{\pi_1^{(i)}(0)}{\pi_1^{(i)}(1)} + \exp\left(-\frac{2\lambda^2}{t\gamma^2}\right).
\end{align*}
Let $x=t \gamma^2/4$ and $y=t \delta^2/2$. Setting $\lambda = -x + \sqrt{2xy + x^2}$ equalizes the exponential terms:
\begin{align*}
\mathbb{E}[\pi^{(i)}_{t+1}(0) \mid Q_1^{(i)}=1] & \leq \left(1+\frac{\pi_1^{(i)}(0)}{\pi_1^{(i)}(1)}\right)\exp\left(-x-y+\sqrt{x^2+2xy}\right)\\&\leq \frac{1}{\pi_1^{(i)}(1)}\exp\Big(-\frac{y^2}{2(x+y)}\Big).
\end{align*}
The second inequality is given by the convex inequality 
$\sqrt{a} - \sqrt{a+b} \leq -\frac{b}{2\sqrt{a+b}}$, for $a=x^2+2xy$ and $b=y^2$.
From the definitions of $x$ and $y$, this yields
\begin{equation*}
\mathbb{E}[\pi^{(i)}_{t+1}(0) \mid Q_1^{(i)}=1]  \leq \frac{1}{\pi_1^{(i)}(1)}\exp\left(-\frac{t \delta^4}{2\gamma^2 + 4 \delta^2}\right). 
\end{equation*}
We conclude by noting that $\pi_t^{(i)}(q'^{(i)}) \geq \pi_t(q')$.
\end{proof}

\section{Dynamical Environment}
\label{sec:dynamic}

We now model a situation where the quality $Q$ may change over time.
We  consider a general Markovian model given by the transition matrix $P$. Moreover, at each time step, the quality might change with probability at most $\eta\in(0,1)$:
\begin{equation}
\label{eq:dynamic}
\begin{gathered}
\mathbb{P}\left( Q_{t+1} = q' \mid Q_t = q \right) = P_{q,q'}, \\
\text{with } P(q,q) \geq 1-\eta \text{ for all } q\in\mathcal{Q}.
\end{gathered}
\end{equation}
The use of a Markovian model is rather usual in such dynamical models and is reasonable when quality does not slowly drift, but rather  changes abruptly; this can happen, for example, after a change of owner, employees, contractor of an establishment or after a technological advance in case of high tech products.

Assuming that the diagonal terms of the transition matrix $P$ are large ensures that changes of quality are rare. 
Consumers thus have some time to learn the current quality of the product.

\medskip

Studying the convergence of the posterior is irrelevant, as the quality regularly changes. Instead, we measure the quality of the posterior variations in term of the total utility loss 
\begin{equation}
\label{eq:regret}
\regret\coloneqq\sum_{t=1}^T \mathbb{E}[ r(Q_t, \theta_t)_+ - u_t],
\end{equation}
also known as ``regret''\footnote{Alternatively, it is the sum of the individual regret of all consumers.}.The first term $r(Q_t, \theta_t)_+$  corresponds to the utility a consumer would get if she knew the quality $Q_t$, whereas $u_t$ is the utility she actually gets. 

\begin{lemma}
\label{lemma:lossequivdiscrete}
If $r$ is $\diffbound$-Lipschitz in its first argument for any $\theta\in \Theta$, \ie 
$|r(q, \theta)-r(q', \theta)| \leq \diffbound \| q-q' \|_1$
for any $q,q' \in \mathcal{Q}$, we have
\begin{equation*}
\regret \leq \diffbound\sum_{i=1}^d\sum_{t=1}^T \mathbb{E}[1-\pi_t^{(i)}(Q_t^{(i)})].
\end{equation*}
\end{lemma}

\cref{lemma:lossequivdiscrete}, proved in \cref{app:proofs_dynamic}, shows that bounding the cumulated estimation error $\sum_{t=1}^T \mathbb{E}[1-\pi_t^{(i)}(Q_t^{(i)})]$ for each coordinate is sufficient to bound the total regret. 

\medskip

In this whole section, we use the notation  $g(T, \eta) = \bigO{f(T, \eta)}$ if there exists a positive constant $c$ such that for all~${T \in \mathbb{N}^*}$ and $\eta \in (0,1)$, $g(T, \eta) \leq c f(T, \eta)$. 
Moreover, we write $g(T, \eta) = \Omega(f(T, \eta))$ if $f(T, \eta) = \mathcal{O}(g(T, \eta))$.

\medskip


We consider in this section consumers who have perfect knowledge of the model, \ie they know that the quality might change following \eqref{eq:dynamic}.
Recall that the prior is assumed uniform on $\mathcal{Q}$.
If $G$ is defined as in \eqref{eq:defnG}, the posterior update is  given  by
\begin{equation}\label{eq:posteriorupdatedyn}
\pi_{t+1}(q) = \sum_{q' \in \mathcal{Q}} P(q,q')\frac{G\left(Z_t,\pi_t, q'\right)}{G\left(Z_t,\pi_t\right)}\pi_t(q').
\end{equation}
The effect of the old reviews is mitigated by the multiplications by the transition matrix $P$. Consumers thus value more recent reviews in this model, as wished in its design.
By induction, the previous equality leads to the following expression.
\begin{equation}\label{eq:posteriordynamic}
\pi_{t+1}(q) = \sum_{\substack{(q_s) \in \mathcal{Q}^t \\ q_{t+1} = q}} \pi_1(q_1)\prod_{s=1}^{t} P(q_s, q_{s+1}) \frac{G(Z_s, \pi_s, q_s)}{G(Z_s, q_s)}.
\end{equation}
This expression is more complex than the one in the stationary case, leading to a more intricate proof of error bounds. 
The multidimensional setting becomes significantly intricate in the dynamic case, as the transition matrix entails correlations between the different dimensions. To overcome this difficulty, we first bound the estimation error for a simpler, imperfect Bayesian estimator ignoring these correlations. It directly leads to a bound on the true utility loss, by optimality of the Bayesian estimator.

\cref{thm:dynamic-binary} below shows that the cumulated loss is of order $\ln(2/\eta)\eta T$. 
Perfect learners, who could directly observe $Q_{t-1}$ before making the decision at time $t$, would still suffer a loss of order $\eta T$ as there is a constant uncertainty $\eta$ about the next step quality\footnote{Equivalently, there are $\eta T$ changes of quality and a positive loss is unavoidable after each change.}. 
\cref{thm:dynamic-binary} thus shows that the cost of learning is just a logarithmic factor in the dynamical setting.%

\begin{theorem}
\label{thm:dynamic-binary} 
If $r$ is $\diffbound$-Lipschitz, then $\regret = \bigO{Md\ln\left(2/\eta\right)\eta T}$.
Moreover, if $\eta T = \Omega(1)$, there is some $\diffbound$-Lipschitz reward $r$ and some transition matrix $P$ verifying the conditions of \cref{eq:dynamic} such that $\regret = \Omega(Md\ln(2/\eta)\eta T)$.
\end{theorem}

The hidden constants in the $\bigO{\cdot}$ and $\Omega(\cdot)$ above only depend on the values of $\delta$ and $\gamma$ defined in \cref{thm:stationary_discrete}.

The proof of \cref{thm:dynamic-binary} is divided into two parts: first, the upper bound $\regret = \bigO{Md\ln(2/\eta)\eta T }$ and, second, the lower bound $\regret=\Omega\left(Md\ln(2/\eta)\eta T  \right)$. The proof of the lower bound is postponed to \cref{app:proofs_dynamic_binary}. 

\medskip

The assumption $\eta T = \Omega(1)$ guarantees that changes of quality actually have a non-negligible chance to happen in the considered time window. Without it, we  would be back to  the stationary case. 
In the extreme case $\eta T \approx 1$, the error is thus of order $\ln(T)$ against $1$ in the stationary setting. 
This larger loss is actually the time needed to achieve the same precision in posterior belief anew after a change of quality. 
Indeed, let the posterior be very close to the true quality $q$, \ie $\pi_t(q') \approx 0$ for $q'\neq q$;  if the quality suddenly changes to $q'$, it will take a while to have a correct estimation again, \ie  to get $\pi_t(q') \approx 1$. %

\medskip

It can directly be seen from the proof that a similar regret bound also holds for imperfect consumers, who do not know the exact transition matrix $P$, but only an upper bound of the constant $\eta$ defined in \cref{eq:dynamic}. Similarly, the results also hold for dynamic transition matrices $P_t$ as long as $P_t(q,q) \geq 1-\eta$ for any $q$ and $t$. This describes a more general setting where a monopolist influences (with limits) the changes of quality that happen.

\subsection{Proof of the Upper Bound.}

In order to prove that $\regret = \bigO{Md\ln(2/\eta)\eta T}$,
we actually show the result marginally on each dimension, \ie for any $i \in [d]$
\begin{equation}
\label{eq:marginalupperbound}
\sum_{t=1}^T 1-\pi^{(i)}_t(q^{(i)}) = \bigO{\ln(2/\eta)\eta T}.
\end{equation}
\cref{lemma:lossequivdiscrete} then directly leads to the upper bound.
To prove \cref{eq:marginalupperbound}, we first consider another $\mathcal{H}_t$-measurable estimator defined for any $i$ by
\begin{equation}
\label{eq:tildeestim}
\tilde{\pi}_1^{(i)} = \pi_1^{(i)} \quad \text{ and } \quad \tilde{\pi}_{t+1}^{(i)}(q^{(i)}) = (1-2\eta) \frac{G^{(i)}(Z_t^{(i)}, \pi_t, (q^{(i)}))}{G^{(i)}(Z_t^{(i)}, \pi_t)}  \tilde{\pi}_{t}^{(i)}(q^{(i)}) + \eta.
\end{equation}
The estimator $\tilde{\pi}_t$ can be seen as the Bayesian estimator, for the worst case of transition matrix, where each feature $i$ changes with probability $\eta$ at each step. 
As perfect Bayesian consumers' decisions minimize the utility loss among the classes of $\mathcal{H}_t$ measurable decisions, having an $\bigO{\ln(2/\eta)\eta T}$ error for $\tilde{\pi}_t^{(i)}$ directly yields \cref{eq:marginalupperbound}.


\medskip

To prove \cref{eq:marginalupperbound}, we partition $\mathbb{N}^*$  into blocks $[t_k^{(i)}+1, t_{k+1}^{(i)}]$ of fixed quality (for the $i$-th coordinate) and show that the error of  $\tilde{\pi}_t^{(i)}$ on each block individually is $\bigO{\ln(2/\eta)}$:
\begin{equation}\label{eq:defntn}
t_1^{(i)} \coloneqq 0 \quad \text{and} \quad t_{k+1}^{(i)} \coloneqq \min\left\{ t > t_k^{(i)} \mid Q_{t+1}^{(i)} \neq Q_{t_k^{(i)}+1}^{(i)} \right\}. 
\end{equation}
%
Define the stopping time \begin{equation}\label{eq:defntau}
\tau_k^{(i)} \coloneqq \min \bigg(\Big\{ t \in [t_k^{(i)}+1, t_{k+1}^{(i)}] \ \Big| \ \frac{\tilde\pi_t^{(i)}(1)}{\tilde\pi_t^{(i)}(0)} \geq 1\Big\} \cup \{t_{k+1}^{(i)}\}\bigg).
\end{equation}
This is the first time\footnote{It is set as the largest element of the block if such a criterion is never satisfied.} in block $k$ where the posterior belief of the true quality (for $\tilde\pi^{(i)}_t$) exceeds the one of the wrong quality. The error on the block is then decomposed as the terms before~$\tau_k^{(i)}$, which contribute to at most $1$ per timestep, and the terms after $\tau_k^{(i)}$. The remaining of the proof, postponed to \cref{app:proofupperbound} due to space constraints, bounds the former term using \cref{lemma:boundtau} below, while \cref{lemma:evolaftertau} allows to bound the latter. The proofs of these lemmas are respectively given in \cref{proof:boundtau,proof:evolaftertau}.
\begin{lemma}\label{lemma:boundtau}
For any $k$, 
\begin{equation*}
\mathbb{P}\left[\tau_k^{(i)}-t_k^{(i)} \geq 2+\frac{2\gamma^2+4\delta^2}{\delta^4}\ln\parens*{\frac{1}{\eta}}\right] \leq \eta,
\end{equation*} 
where $\delta$ and $\gamma$ are defined as in \cref{thm:stationary_discrete}.
\end{lemma}
%
%
\begin{lemma}\label{lemma:evolaftertau}For any $k \in \mathbb{N}^*$ and $t \in [\tau_k^{(i)}, t_{k+1}^{(i)} ] $,
\begin{equation*}
\mathbb{E}\left[\frac{1}{\tilde\pi^{(i)}_t(Q_{t}^{(i)})}  \ \Big\vert \ \tau_k^{(i)}, (t_n^{(i)})_n \right] \leq 2.
\end{equation*}
\end{lemma}


\subsection{Proof of the Lower Bound.} The proof of the lower bound is postponed to \cref{app:proofs_dynamic_binary}. The idea is  that the posterior cannot converge faster than exponentially on a single block. Thus, if the posterior converged in the last block, \eg $\pi_t(q') \approx \eta$ in a block of quality $q$, then it would require a time $\ln(2/\eta)$ before $\pi_t(q') \geq 1/2$ in the new block of quality $q'$, leading to a loss at least $\ln(2/\eta)$ on this block.
\section{Naive Learners} 
\label{sec:imperfect_learners}

In \cref{sec:dynamic} we showed that learning occurs for Bayesian consumers who are perfectly aware of the environment, and especially of its dynamical aspect. 
In some learning problems, Bayesian learners can still have small regret, despite having an imperfect knowledge of the problem parameters or even ignoring some aspects of the problem.

This section shows that awareness of the problem's dynamical structure is essential here. In particular, naive learners incur a considerable utility loss.

In the following, we consider the setting described in \cref{sec:dynamic} with naive learners, \ie  consumers who are unaware of possible quality changes over time. 
As a consequence, their posterior distribution $\pi_t^{\imp}$ follows the exact same update rule as in the stationary case: 
\begin{equation*}
\pi_{t+1}^\imp(q) = \frac{G(Z_t, \pi_t^\imp, q)}{G(Z_t, \pi_t^\imp)}\pi_t^\imp(q).
\end{equation*} 
The regret for naive learners is then
\begin{gather*}
\sum_{t=1}^T\mathbb{E}\left[ r(Q_t, \theta_t)_+ - u_t^\imp \right],\\
\text{where }u_t^{\imp} = r(Q_t, \theta_t)\one_{\sum_{q\in\mathcal{Q}}\pi_t^\imp(q)r(q,\theta_t)\geq0}.
\end{gather*}
$u_t^{\imp}$ is the utility achieved by naive learners who make their decisions based on $\pi_t^\imp$.

\cref{thm:imperfectloss} below states that the utility loss for naive learners is non-negligible, \ie of order $T$, which displays the significance of taking into account the dynamical structure of the problem in the learning process.
\begin{theorem}\label{thm:imperfectloss}
If $\eta T = \Omega(1)$, then there is some $\diffbound$-Lipschitz reward $r$ and some transition matrix $P$ verifying the conditions given by \cref{eq:dynamic} such that 
\begin{equation*}
\regret^\imp = \Omega(MdT).
\end{equation*}
\end{theorem}
The proof of \cref{thm:imperfectloss} can be found in \cref{app:proofs-dynamic-imperfect} and bears similarities with the proof of the lower bound in \cref{thm:dynamic-binary}. 
The posterior of naive learners converges quickly to the true quality on a single block. 
Because of this, after a change of quality, it takes a long time before the posterior belief of naive learners becomes accurate again with respect to the new quality.

\section{Conclusions and future work}
\label{sec:conclusion}

This work was a first attempt to use a change-point framework for social learning in online markets with reviews when the qualities of a product vary over time. 
Leveraging convergence rates in the stationary setting, we obtained regret bounds in a dynamical model when the quality space is $\{0,1\}^{d}$; our bound is tight.
For more general quality spaces, determining the incurred regret is a much harder problem. 
We propose some partial results for continuous qualities in \cref{app:continuous}, where we provide lower and upper bounds of the regret. 
Unfortunately the gap between these bounds remains to be closed.

Many other directions also remain open for review based markets. 
For instance, it would be interesting to study a model with a slowly drifting quality, rather than abrupt changes. 

In this work we only focus on the consumer side, but the seller can also adaptively set the price (or the quality) of the item. 
What is a good seller strategy in this case?
The selling platform can also design the format of the feedback given by the consumers. Determining the format which allows the best possible convergence rate might also be of great interest in practice. This work thus developed the techniques to analyze the customers' behavior in social learning in non-stationary environments and thus opens possibilities for future developments on pricing or platform design to improve the revenue of the platform and/or the experience of the customer.

Considering perfect Bayesian consumers might be unrealistic. 
In reality, consumers have limited computation capacity or can be risk averse, leading to different behaviors. 
Studying the effects of these limitations is also of great interest.

\paragraph{Acknowledgements.}
This research project received partial support from the COST action GAMENET.
Marco Scarsini is a member of INdAM-GNAMPA.
His work was partially supported by the Italian MIUR PRIN 2017 Project ALGADIMAR ``Algorithms, Games, and Digital Markets.'' 
He gratefully acknowledges the kind hospitality of Ecole Normale Sup\'erieure Paris-Saclay, where this project started. 
This work was completed while Etienne Boursier was a PhD student at Centre Borelli, ENS Paris-Saclay, France.
\printbibliography
\appendix
\newpage
\section{List of symbols}
\label{se:symbols}
This is a list of the symbols used throughout the paper.
\begin{longtable}{p{.25\textwidth} p{.7\textwidth}}

$f$ & feedback function \\
$G(z, \pi, q)$ & $\mathbb{P}[Z_t=z \mid \pi_t=\pi, Q_t=q]$, defined in 
\eqref{eq:defnG} \\
$G(z, \pi)$ & $\mathbb{E}_{q \sim \pi}[G(z, \pi, q)]$, defined in \eqref{eq:defnGzpi} \\
$G^{(i)}(z^{(i)}, \pi, q^{(i)})$ & $\mathbb{P}[Z_{t}^{(i)}=z^{(i)} \mid \pi_t=\pi, Q_t^{(i)}=q^{(i)}]$, defined in \eqref{eq:defnGi} \\
$G^{(i)}(z^{(i)}, \pi)$ & $\mathbb{E}_{q \sim \pi}[G^{(i)}(z^{(i)}, \pi, q^{(i)})]$, defined in \eqref{eq:defnGizpi} \\
$h(x)$ & $\parens*{\frac{\eta}{2} + \frac{1 -\eta}{x}}^{-1}$, defined in \eqref{eq:h} \\
$\mathcal{H}_{1}$ & $\varnothing$ \\
$\mathcal{H}_{t}$ & history at time $t$ \\
$\widetilde{\mathcal{H}}_{t}$ & reduced history at time $t$ \\
$\kl{\mu,\nu}$ & Kullback-Leibler divergence, defined in \eqref{eq:KL} \\
$\diffbound$ & Lipschitz constant of $r$, defined in \cref{lemma:lossequivdiscrete} \\
$p^{(i)}$ & price of product $i$ \\
$Q_{t}$ & quality at time $t$ \\
$r$ & reward function \\
$\regret$ & $\sum_{t=1}^T \mathbb{E}[ r(Q_t, \theta_t)_+ - u_t]$, regret, defined in \eqref{eq:regret} \\
$u_{t}$ & utility of consumer $t$ \\
$x^{(i)}$ & $i$-th component of the vector $x$ \\
$X_t$ & $\ln\left(\frac{G^{(i)}(Z_t^{(i)}, \pi_t, 1)}{G^{(i)}(Z_t^{(i)}, \pi_t, 0)} \right) - \kl{G^{(i)}(\cdot, \pi_t, 1), G^{(i)}(\cdot, \pi_t, 0)}$, defined in the proof of \cref{thm:stationary_discrete} \\
$Z_{t}$ & review of consumer $t$ \\

$\gamma$ & $2\max_{\substack{i\in [d], \pi \in \mathcal{P}(\mathcal{Q}), z \in \mathcal{Z}}}  \left| \ln\left(\frac{G^{(i)}(z^{(i)}, \pi, 1)}{G^{(i)}(z^{(i)}, \pi, 0)}\right) \right|$, defined in \eqref{eq:gamma-i} \\
$\delta$ & $\min_{\substack{i\in [d], \pi \in \mathcal{P}(\mathcal{Q})}} \sum_{z \in \mathcal{Z}}| G^{(i)}(z^{(i)}, \pi, 1) - G^{(i)}(z^{(i)}, \pi, 0)|
$, defined in \eqref{eq:delta-i} \\

$\varepsilon_t$ & noise of consumer $t$ \\
$\eta$ & probability that the quality changes, defined in \eqref{eq:dynamic} \\ 
$\theta_{t}$ & preference of consumer $t$ \\
$\xi_{t}$ & binary vector indicating which features are reviewed by consumer~$t$ \\
$\pi_{1}$ & prior distribution of the quality vector \\
$\pi_{t}$ & posterior distribution of the quality vector at time $t$ \\
\end{longtable}
\section{Additional proofs}
\label{se:additional-proofs}
This section contains detailed proofs of lemmas and theorems postponed to the Appendix.

\subsection{Proof of Lemma~\ref{lemma:lossequivdiscrete}}\label{app:proofs_dynamic}

The inequality actually holds for any term of the sum when conditioned on $\pi_t$, \ie $\mathbb{E}[r(Q_t, \theta_t)_+ - u_t \mid \pi_t] \leq M \sum_{i=1}^d(1-\pi_t^{(i)}(Q_t^{(i)}))$, which directly implies \cref{lemma:lossequivdiscrete}.
By definition, $u_t = r(Q_t, \theta_t)\one_{\sum_{q\in\mathcal{Q}}\pi_t(q)r(q,\theta_t)\geq0}$ and so it comes
\begin{align*}
r(Q_t, \theta_t)_+ - u_t & = r(Q_t, \theta_t) \left(\one_{r(Q_t, \theta_t) \geq0} -\one_{\sum_{q\in\mathcal{Q}}\pi_t(q)r(q,\theta_t)\geq0}\right) \\
& = r(Q_t, \theta_t) \left(\one_{r(Q_t, \theta_t) \geq0 \geq\sum_{q\in\mathcal{Q}}\pi_t(q)r(q,\theta_t) } -\one_{\sum_{q\in\mathcal{Q}}\pi_t(q)r(q,\theta_t)\geq0\geq r(Q_t, \theta_t)}\right)\\
& \leq \left| r(Q_t, \theta_t) -\sum_{q\in\mathcal{Q}}\pi_t(q)r(q,\theta_t)\right| \\
& = \left| \sum_{q\in\mathcal{Q}}\pi_t(q)(r(Q_t, \theta_t)-r(q,\theta_t))\right|\\
& \leq \sum_{q\in \mathcal{Q}} \pi_t(q)\left| r(Q_t,\theta_t)- r(q, \theta_t)\right| \\
& \leq M \sum_{q\in \mathcal{Q}} \pi_t(q) \|Q_t - q\|_1 \\
& = M \sum_{i=1}^d (1-\pi_t^{(i)}(Q_t^{(i)})).
\end{align*}

\subsection{Proof of the upper bound in \cref{thm:dynamic-binary}}
\label{app:proofupperbound}
Recall that we partitioned $\mathbb{N}^*$ into blocks $[t_k^{(i)}+1, t_{k+1}^{(i)}]$ of fixed quality and we also defined in \cref{eq:defntau} the stopping time $\tau_k^{(i)}$ as the first time of the block where the posterior belief of the quality (for $\tilde{\pi}_t^{(i)}$) becomes larger than $\frac{1}{2}$.

We now only aim at bounding the estimation error on a single block $k$ and we assume w.l.o.g.\ that $Q_t^{(i)}=1$ on this block.

Similarly to the proof of \cref{thm:stationary_discrete}, Azuma-Hoeffding's inequality on a single block $k$ leads to
\begin{equation}\label{eq:evolaftertau1}
\mathbb{E}\left[\prod_{s=n}^{t-1} \frac{G^{(i)}(Z_s^{(i)}, \pi_s, 0)}{G^{(i)}(Z_s^{(i)}, \pi_s, 1)} \,\Big|\, \pi_{n}, \forall s \in [n, t-1], Q_s^{(i)} = 1\right] \leq \exp\left(-\frac{(t-n) \delta^4}{2\gamma^2+4\delta^2}\right).
\end{equation}%
Also, note that \cref{eq:tildeestim} leads to 
\begin{equation*}
\frac{G^{(i)}(Z_t^{(i)}, \pi_t, 1)}{G^{(i)}(Z_t^{(i)}, \pi_t)} \leq \frac{\tilde\pi^{(i)}_{t+1}(1)}{(1-2\eta)\tilde\pi^{(i)}_t(1)}.
\end{equation*}
By induction, we get
\begin{equation}
\label{eq:prod-induc}
\prod_{s=n}^{t-1}\frac{G^{(i)}(Z_s^{(i)}, \pi_s, 1)}{G^{(i)}(Z_s^{(i)}, \pi_s)} \leq \frac{1}{\tilde\pi^{(i)}_n(1)(1-2\eta)^{t-n}}.
\end{equation}
Multiplying the left hand side of \eqref{eq:evolaftertau1} by the left hand side of \eqref{eq:prod-induc}, we obtain
\begin{equation}\label{eq:prodG1}
\mathbb{E}\left[\prod_{s=n}^{t-1} \frac{G^{(i)}(Z_s^{(i)}, \pi_s, 0)}{G^{(i)}(Z_s^{(i)}, \pi_s)} \ \Big\vert\ \pi_{n}, \forall s \in [n, t-1], Q_s^{(i)} = 1\right] \leq \frac{(1-2\eta)^{-(t-n)}}{\tilde\pi_n^{(i)}(1)}\exp\left(-\frac{(t-n) \delta^4}{2\gamma^2+4\delta^2}\right).
\end{equation}
%
Similarly to \cref{eq:posteriordynamic},  starting from $n_0\geq 1$, for the $i$-th coordinate it can be shown that
\begin{align*}
\label{eq:posteriordyn2}
\tilde\pi_{t+1}^{(i)}(q^{(i)}) &= (1-2\eta)^{t-n_0+1} \tilde\pi^{(i)}_{n_0}(q^{(i)}) \prod_{s=n_0}^t \frac{G^{(i)}\left(Z_s^{(i)},\pi_s, q^{(i)}\right)}{G^{(i)}\left(Z_s^{(i)},\pi_s\right)} \\
&\quad+ \eta \sum_{s=0}^{t-n_0}(1-2\eta)^{s}\!\!\!\!\prod_{l=t-s+1}^{t} \frac{G^{(i)}(Z_l^{(i)}, \pi_l, q^{(i)})}{G^{(i)}(Z_l^{(i)}, \pi_l)}.
\end{align*}%
Define $A_{\tau_k^{(i)}}^t \coloneqq \left\{\forall s \in [\tau_k^{(i)}, \tau_k^{(i)}+t], Q_s^{(i)} = 1 \right\}$. Combining this formula with \cref{eq:prodG1}, we obtain
\begin{align*}
\mathbb{E}\bracks*{\tilde\pi_{\tau_k^{(i)}+t}^{(i)}(0) \mid \mathcal{H}_{\tau_k^{(i)}}, A_{\tau_k^{(i)}}^t}\leq& \frac{\tilde\pi_{\tau_k^{(i)}}^{(i)}(0)}{\tilde\pi^{(i)}_{\tau_k^{(i)}}(1)}\exp\left(-\frac{t \delta^4}{2\gamma^2+4\delta^2}\right)\\ 
&+ 2\eta \sum_{s=0}^{t-1} \mathbb{E}\left[\frac{1}{\tilde\pi^{(i)}_{\tau_k^{(i)}+t-s}(1)} \,\Big|\, \mathcal{H}_{\tau_k^{(i)}},  A_{\tau_k^{(i)}}^t\right]\exp\left(-\frac{s \delta^4}{2\gamma^2+4\delta^2}\right).
\end{align*}
Thanks to \cref{lemma:evolaftertau}, 
\begin{equation*}
\mathbb{E}\left[\frac{1}{\tilde\pi^{(i)}_{\tau_k^{(i)}+t-s}(1)}\,\Big|\, \mathcal{H}_{\tau_k^{(i)}}, A_{\tau_k^{(i)}}^t\right] \leq 2 \quad\text{and}\quad \frac{\tilde\pi^{(i)}_{\tau_k^{(i)}}(0)}{\tilde\pi^{(i)}_{\tau_k^{(i)}}(1)} \leq 1,
\end{equation*} 
so that
\begin{equation}\label{eq:dynamicconcentr}
\begin{split}
\mathbb{E}\Big[\tilde\pi^{(i)}_{\tau_k^{(i)}+t}(0) \mid \mathcal{H}_{\tau_k^{(i)}},  A_{\tau_k^{(i)}}^t\Big] & \leq \exp\left(-\frac{t \delta^4}{2\gamma^2+4\delta^2}\right) +  4\eta\sum_{s=0}^{t-1} \exp\left(-\frac{s \delta^4}{2\gamma^2+4\delta^2}\right) \\
&\leq \exp\left(-\frac{t \delta^4}{2\gamma^2+4\delta^2}\right) + \frac{4\eta}{1-\exp\left(-\frac{\delta^4}{2\gamma^2+4\delta^2}\right)}.
\end{split}
\end{equation}
Finally, the estimation error for $\tilde\pi_t^{(i)}$ incurred during the block $k$ is at most 
\begin{equation*}
\tau_k^{(i)} - t_k^{(i)} + \sum_{t=0}^{t_{k+1}^{(i)}-t_k^{(i)}-1} \parens*{\exp\left(-\frac{t \delta^4}{2\gamma^2+4\delta^2}\right) + \frac{4\eta}{1-\exp\left(-\frac{\delta^4}{2\gamma^2+4\delta^2}\right)}},
\end{equation*}
\ie it is of order $\tau_k^{(i)} -t_k^{(i)} +  \eta (t_{k+1}^{(i)}-t_k^{(i)})$.
\cref{lemma:boundtau} then yields 
\begin{equation*}
\mathbb{E}[\tau_k^{(i)} -t_k^{(i)} \mid (t_n)_n] \leq 2+\frac{2\gamma^2+4\delta^2}{\delta^4}\ln\parens*{\frac{1}{\eta}} +  \eta(t_{k+1}^{(i)}-t_k^{(i)}).
\end{equation*}
Thus in expectation, given $(t_n)_n$, the estimation error of $Q_t^{(i)}$ over the block $k$ for $\tilde\pi_t$ is of order
$
\ln(2/\eta) + \eta(t_{k+1}^{(i)}-t_k^{(i)})$.
Note  that $t_{k+1}^{(i)}-t_k^{(i)}$ is stochastically dominated by a geometric distribution of parameter $\eta$. In expectation the number of blocks counted before $T$ is thus $\bigO{\eta T}$ and summing over all these blocks yields 
\begin{equation*}
\sum_{t=1}^T \mathbb{E}[1-\tilde\pi_t^{(i)}(Q_t^{(i)})] = \bigO{\ln(2/\eta)\eta T}.
\end{equation*}
When summing over all coordinates, this implies that the regret incurred by the estimator $\tilde\pi_t$ is of order $ \bigO{Md\ln(2/\eta)\eta T}$. Since the exact estimator $\pi_t$ minimizes the expected utility loss among the class of all $\mathcal{H}_t$-measurable estimators, the upper bound of \cref{thm:dynamic-binary} follows.

\subsection{Proof of \cref{lemma:boundtau}}
\label{proof:boundtau}
As a consequence of the posterior update of $\tilde\pi_t$ given by \cref{eq:tildeestim}, for $t+1\leq \tau_k^{(i)}$,
\begin{equation*}
\tilde\pi_{t+1}^{(i)}(0) \leq  \frac{G^{(i)}(Z_t^{(i)}, \pi_t, 0)}{G(Z_t, \pi_t)}\tilde\pi_t(0) \quad \text{ and } \quad \tilde\pi_{t+1}^{(i)}(1) \geq  \frac{G^{(i)}(Z_t^{(i)}, \pi_t, 1)}{G(Z_t, \pi_t)}\tilde\pi_t(1).
\end{equation*}
We then get by induction
\begin{equation}
\label{eq:linearupdate}
\frac{\tilde\pi^{(i)}_{t+1}(0)}{\tilde\pi^{(i)}_{t+1}(1)} \leq \frac{1}{\eta} \prod_{s=t_k^{(i)}+1}^{t}\frac{G^{(i)}(Z_s^{(i)}, \pi_s, 0)}{G^{(i)}(Z_s^{(i)}, \pi_s, 1)},
\end{equation}
as $\tilde\pi^{(i)}_{t_k^{(i)}+1}(1) \geq  \eta$.
For $n = \left\lceil \frac{2\gamma^2+4\delta^2}{\delta^4}\ln\parens*{\frac{1}{\eta}} \right\rceil$, it has been shown in the proof of \cref{thm:stationary_discrete} that:
\begin{equation*}
\mathbb{P}\bracks*{\prod_{s=t_k^{(i)}+1}^{t_k^{(i)}+n}\frac{G(Z_s^{(i)}, \pi_s, 0)}{G(Z_s^{(i)}, \pi_s, 1)} > \eta  \ \Big\vert \ \pi_{t_k^{(i)}+1}, \forall s \in [t_k^{(i)}+1, t_k^{(i)}+n], Q_s^{(i)} = 1} \leq \eta.
\end{equation*} 
Note that by definition of $\tau_k^{(i)}$, 
$
\frac{\tilde\pi^{(i)}_{\tau_k^{(i)}}(0)}{\tilde\pi^{(i)}_{\tau_k^{(i)}}(1)} \leq 1$.
The above  concentration inequality and \eqref{eq:linearupdate} imply that $\mathbb{P}[\tau_k^{(i)} - t_k^{(i)} \geq n+1] \leq \eta$. 

\subsection{Proof of \cref{lemma:evolaftertau}}
\label{proof:evolaftertau}
By definition of $G^{(i)}$ and the posterior update,  given by \cref{eq:defnGi,eq:tildeestim} respectively, we have
\begin{gather}
\mathbb{E}\bracks*{\frac{1}{\tilde\pi_{t+1}^{(i)}(1)} \ \Big\vert \ Q_t^{(i)} = 1, \mathcal{H}_t} = \sum_{z^{(i)}: z \in \mathcal{Z}} G^{(i)}(z^{(i)}, \pi_t^{(i)}, 1) \, h\left(\frac{G^{(i)}(z^{(i)}, \pi_t)}{G^{(i)}(z^{(i)}, \pi_t,1) \tilde\pi^{(i)}_t(1)} \right),
\nonumber\\
\label{eq:h}
\text{with } h(x) = \frac{1}{\eta + \frac{1 -2\eta}{x}}.
\end{gather}
Note that $h$ is concave on $\mathbb{R}_+^*$, so by Jensen's inequality:
\begin{equation}\label{eq:evoltau1}
\mathbb{E}\bracks*{\frac{1}{\tilde\pi^{(i)}_{t+1}(1)}  \ \Big\vert \ Q_t^{(i)} = 1, \mathcal{H}_t} \leq h\left(\frac{1}{\tilde\pi^{(i)}_t(1)}\right).
\end{equation}

\cref{lemma:evolaftertau} then follows by induction
\begin{align*}
&\mathbb{E}\left[\frac{1}{\tilde\pi_{t+\tau_k^{(i)}+1}^{(i)}(1)} \ \Big\vert \ \tau_k^{(i)}, \forall s \in [\tau_k^{(i)}, t+\tau_k^{(i)}], Q_s^{(i)} = 1 \right] \\
&\qquad \leq \mathbb{E}\left[h\left(\frac{1}{\tilde\pi_{t+\tau_k^{(i)}}^{(i)}(1)}\right) \ \Big\vert \ \tau_k^{(i)}, \forall s \in [\tau_k^{(i)}, t+\tau_k^{(i)}], Q_s^{(i)} = 1 \right]\\
&\qquad \leq  h\left(\mathbb{E}\left[\frac{1}{\tilde\pi^{(i)}_{t+\tau_k^{(i)}}(1)}  \ \Big\vert \ \tau_k^{(i)}, \forall s \in [\tau_k^{(i)}, t+\tau_k^{(i)}], Q_s^{(i)} = 1 \right]\right) \\
&\qquad \leq h\left(2\right) =2.
\end{align*}
The first inequality is a direct consequence of \cref{eq:evoltau1}, the second is Jensen's inequality again, while the third one is obtained by induction using the fact that $h$ is increasing and $\tilde\pi_{\tau_k^{(i)}}^{(i)}(1)\geq \frac{1}{2}$.

\subsection{Proof of the lower bound in \cref{thm:dynamic-binary}}\label{app:proofs_dynamic_binary}

In this proof we consider the following transition matrix:
\begin{equation*}
P(q, q) = 1-\eta \quad \text{ and } P(q, \pmb{1}-q) = \eta,
\end{equation*}
\ie all the features change simultaneously with probability $\eta$ at each round. We also assume that the prior is only split between the vectors $\pmb{0}$ and $\pmb{1}$, \ie the features are either all $0$ or all $1$. If we take the reward function $r(q, \theta)= M\sum_{i=1}^d q_i + \theta_i$, then the regret scales as
\begin{align}
\regret & = \Omega\left(M\sum_{i=1}^d \sum_{t=1}^T \mathbb{E}[1-\pi_t^{(i)}(Q_t^{(i)})]\right) \nonumber \\
& = \Omega\left(Md \sum_{t=1}^T\mathbb{E}[1-\pi_t(Q_t)]\right).\label{eq:lowerregreteq}
\end{align}

In this model, we thus have the following posterior update
\begin{equation}\label{eq:lowerupdate}
\pi_{t+1}(\pmb{1}) = (1-2\eta) \frac{G(Z_t, \pi_t, \pmb{1})}{G(Z_t, \pi_t)} \pi_t(\pmb{1}) + \eta.
\end{equation}

This proof uses a partitioning in blocks as follows
\begin{equation}\label{eq:defntnglobal}
t_1 \coloneqq 0 \quad \text{and} \quad t_{k+1} \coloneqq \min\left\{ t > t_k \mid Q_{t+1} \neq Q_{t_k+1} \right\}.
\end{equation}

Consider the block $k$ and assume w.l.o.g. that $Q_t=\pmb{1}$ for this block. 
Define the stopping time \begin{equation}\label{eq:defntaulower}
\tau_{k} \coloneqq \min \bigg(\Big\{ t \in [t_{k}+1, t_{k+1}] \ \Big| \ \pi_t(\pmb{1}) \geq \frac{1}{2}\Big\} \cup \{t_{k+1}\}\bigg),
\end{equation}
and similarly for $\tau_{k+1}$ (with $\pmb{0}$).

The estimation error incurred during blocks $k$ and $k+1$ is at least $\left(\tau_k -t_k + \tau_{k+1} - t_{k+1}\right)/2$. 

Given the posterior update, $\pi_{t+1}(\pmb{1}) \leq c\pi_t(\pmb{1})$ where $c = 1+\max_{\pi, z} \frac{G(z, \pi, \pmb{1})}{G(z, \pi, \pmb{0})}$. As a consequence, $\tau_{k+1} -t_{k+1} \geq \min\left(-\frac{\ln(2\pi_{t_{k+1}}(\pmb{0}))}{\ln(c)}, t_{k+2}-t_{k+1}\right)$. Assume in the following that $t_{k+2}-t_{k+1} \geq -\frac{\ln( 2\eta)}{\ln(c)}$, so that we actually have $\tau_{k+1} - t_{k+1} \geq-\frac{\ln(2\pi_{t_{k+1}}(0))}{\ln(c)}$.

We now bound $\ln(\pi_{t_{k+1}}(\pmb{0}))$ in expectation. By concavity of the logarithm,%
\begin{equation*}
\mathbb{E}[\ln(\pi_{t_{k+1}}(\pmb{0})) \ | \ (t_n)_n, \tau_k] \leq \ln\left( \mathbb{E}[\pi_{t_{k+1}}(\pmb{0}) \ | \ (t_n)_n, \tau_k]\right).
\end{equation*}

Note that the estimator $\tilde\pi_t$ in the proof of the upper bound is similar to $\pi_t$ for the transition matrix considered here. Equation~\eqref{eq:dynamicconcentr} then yields
\begin{equation*}
\mathbb{E}[\pi_{t_{k+1}}(\pmb{0}) \ | \ (t_n)_n, \tau_k ] \leq  \exp\left(-\frac{(t_{k+1}-\tau_k)\bar\delta^4}{2\bar\gamma^2+4\bar\delta^2}\right) + \frac{ \eta}{1-\exp\left(-\frac{\bar\delta^4}{2\bar\gamma^2+4\bar\delta^2}\right)},
\end{equation*}
where \begin{equation*}
\bar\delta \coloneqq \min_{\pi \in \mathcal{P}(\mathcal{Q})} \sum_{z \in \mathcal{Z}}| G(z, \pi, \pmb{1}) - G(z, \pi, \pmb{0})| \quad \text{ and }\quad
\bar\gamma \coloneqq 2\max_{\pi \in \mathcal{P}(\mathcal{Q}), z \in \mathcal{Z}}  \left| \ln\left(\frac{G(z, \pi, \pmb{1})}{G(z, \pi, \pmb{0})}\right) \right|.\end{equation*}

And so, with $t_{k+2}-t_{k+1} \geq -\frac{\ln( \eta)}{\ln(c)}$,
\begin{adjustwidth}{-0.7in}{-0.7in}
\begin{align*}
\mathbb{E}[\tau_k - t_k + \tau_{k+1} - t_{k+1} \mid (t_n)_n, \tau_k] & \geq \tau_k - t_k + \Omega\left(-\ln\left(\exp\left(-\frac{(t_{k+1}-\tau_k)\bar\delta^4}{2\bar\gamma^2+4\bar\delta^2}\right) + \frac{ \eta}{1-\exp\left(-\frac{\bar\delta^4}{2\bar\gamma^2+4\bar\delta^2}\right)}\right)_+\right)  \\
&  \geq  \tau_k - t_k  + \Omega\left(\left(-\ln\left(\eta\right)-\frac{1}{\eta}\exp\left(-\frac{(t_{k+1}-\tau_k)\bar\delta^4}{2\bar\gamma^2+4\bar\delta^2}\right)\right)_+\right).
\end{align*}
\end{adjustwidth}
Where we used the convex inequality $-\ln(x+y)\geq -\ln(x)-y/x$.

When looking at the variations of the right hand side with $\tau_k$, it is minimized either when $\tau_k=t_k$ or when the second term is equal to $0$, \ie $t_{k+1}-\tau_k= \Omega(\ln(1/\eta))$. Finally this yields when~$t_{k+2}-t_{k+1} \geq -\frac{\ln(2 \eta)}{\ln(c)}$:
\begin{equation}\label{eq:minlossblock}
\begin{aligned}
\mathbb{E}[\tau_k - t_k + \tau_{k+1} - t_{k+1} \mid (t_n)_n] & \geq \Omega\bigg(\min\Big(& \ln( 1/\eta)-\frac{1}{ \eta}\exp\left(-\frac{(t_{k+1}-t_k)\bar\delta^4}{2\bar\gamma^2+4\bar\delta^2}\right),&\\& 
& \ln(1/\eta)+t_{k+1}-t_k&\Big)\bigg).
\end{aligned}
\end{equation}

\paragraph{Case $\eta T \geq 32$.}

Recall that $t_{k+1}-t_k$ are \iid geometric variables of parameter $\eta$. Lemma~\ref{lemma:concentrationgeom} below provides some concentration bound for the sum of such variables. Its proof is given at the end of the section.
\begin{lemma} \label{lemma:concentrationgeom}
Denote by $Y(n,p)$ the sum of $n$ \iid geometric variables of parameter $p$. We have the following concentration bounds on $Y(n,p)$:
\begin{enumerate}
\item For $k\leq 1$ and $kn/p \in \mathbb{N}$, $\mathbb{P}[Y(n,p)<kn/p] \leq \exp\left(-\frac{(1-1/k)^2kn}{1+1/k}\right)$.
\item For $k \geq 1$ and $kn/p \in \mathbb{N}$, $\mathbb{P}[Y(n,p)>kn/p] \leq \exp\left(-\frac{(1-1/k)^2kn}{2}\right)$.
\end{enumerate}
\end{lemma}

Let $\beta \in [\frac{1}{4}, \frac{1}{2}]$ such that $\beta  \eta T \in 2\mathbb{N}$ and note that $\one_{t_{k+1}-t_k \geq x}$ follows a Bernoulli distribution of parameter smaller than $(1- \eta)^{\lceil x \rceil}$. We then have the following concentration bounds:
\begin{equation}\label{eq:geombound1}
\mathbb{P}\left[ \sum_{k=1}^{\beta  \eta T} {t_{k+1}-t_k  > T} \right] \leq \exp\left(-\frac{(1-\beta)^2  \eta T}{2}\right) \leq \exp\left(-\frac{ \eta T}{8}\right) \leq e^{-4}. \end{equation}
and
\begin{equation} \label{eq:geombound2}
\begin{gathered}
\mathbb{P}\left[ \sum_{k=1}^{\beta  \eta T/2}\!\!\!\! \one_{t_{2k+1}-t_{2k} \geq \frac{1}{ \eta}} \one_{t_{2k+2}-t_{2k+1} \geq -\frac{\ln(2 \eta)}{\ln(c)}} \leq \frac{\beta  \eta T}{4} (1- \eta)^{\frac{1}{ \eta} - \frac{\ln(2 \eta)}{\ln(c)}} \right] \\ \leq \exp\left(-\frac{\beta \eta T (1- \eta)^{\frac{1}{ \eta} - \frac{\ln(2 \eta)}{\ln(c)}}}{16}\right).\end{gathered}
\end{equation}

The first bound is a direct consequence of Lemma~\ref{lemma:concentrationgeom} while the second one is an application of Chernoff bound to Bernoulli variables of parameter $(1- \eta)^{\lceil\frac{1}{ \eta} \rceil - \lfloor \frac{\ln(2 \eta)}{\ln(c)}\rfloor}$.
Recall that we only consider small $ \eta$. We can thus assume that $ \eta$ is small enough so that $\frac{1}{ \eta} \geq -\frac{\ln(2 \eta)}{\ln(c)}$. The second bound then becomes:
\begin{align*}
\mathbb{P}\left[ \sum_{k=1}^{\beta  \eta T/2} \one_{t_{2k+1}-t_{2k} \geq \frac{1}{\eta}} \one_{t_{2k+2}-t_{2k+1} \geq -\frac{\ln(2 \eta)}{\ln(c)}} \leq \frac{\beta  \eta T}{4} (1- \eta)^{\frac{2}{ \eta}} \right] \leq \exp\left(-\frac{\beta  \eta T (1- \eta)^{\frac{2}{ \eta}}}{16}\right).
\end{align*}%

Note that for any $x \in (0,\frac{1}{2})$, $e^{-3} \leq (1-x)^{2/x}$, so that the last inequality implies for $\eta\leq\frac{1}{2}$
\begin{align*}
\mathbb{P}\left[ \sum_{k=1}^{\beta \eta T/2} \one_{t_{2k+1}-t_{2k} \geq \frac{1}{\eta}} \one_{t_{2k+2}-t_{2k+1} \geq -\frac{\ln(\eta)}{\ln(c)}} \leq \frac{\beta \eta T}{4} e^{-3} \right] &\leq \exp\left(-\frac{\beta \eta T e^{-3}}{16}\right) \\& \leq \exp\left(-\frac{7 e^{-3}}{8}\right).
\end{align*}%

Now note that $e^{-4}+e^{-\frac{7 e^{-3}}{8}} < 1$ so that neither the event in Equation~\eqref{eq:geombound1} nor in Equation~\eqref{eq:geombound2} hold with some constant probability. In that case, Equation~\eqref{eq:geombound1} means that the $\beta \eta T$ first blocks fully count in the regret. Equation~\eqref{eq:geombound2} implies that Equation~\eqref{eq:minlossblock} holds for at least $\Omega( \eta T)$ pairs of blocks and for each of them, the incurred error is at least $\Omega(\ln(1/ \eta))$. This finally implies that $\loss = \Omega(\ln(1/ \eta)  \eta T)$ and similarly for the regret. We conclude by summing over all the coordinates.
\paragraph{Case $ \eta T \leq 32$.}
Since $ \eta T = \Omega(1)$, we can consider a constant $c_0 >0$ such that $ \eta T > c_0$. In that case, the desired bound can actually be obtained on the two first blocks only. Assume w.l.o.g. for simplicity that $T$ is a multiple of $4$.
\begin{align*}
\mathbb{P}\left( t_1-t_0  \in [T/4, T/2] \text{ and } t_2  - t_1  \in [T/4, T/2]\right) & = \left( (1- \eta)^{T/4} - (1- \eta)^{T/2} \right)^2 \\
& = e^{\frac{T}{2}\ln(1- \eta)}(1-e^{\frac{T}{4}\ln(1- \eta)})^2\\
& = e^{-\frac{\eta T}{2}}(1-e^{-\eta T/2})^2.
\end{align*}%
With a positive probability depending only on $c_0$, the two first blocks are completed before $T$, $t_1  - t_0  \geq T/4$ and $t_2 -t_1  \geq T/4$. 
Assuming w.l.o.g. that $ \eta$ is small enough so that $T/4 \geq -\frac{\ln(2 \eta)}{\ln(c)}$, Equation~\eqref{eq:minlossblock} then gives that the loss incurred during the two first blocks is $\Omega\left( \ln(1/ \eta)\right)$. As $ \eta T = \bigO{1}$ in this specific case, this still leads to $$\sum_{t=1}^T \mathbb{E}[1-\pi_t(Q_t)] = \Omega\left( \ln(1/ \eta) \eta T\right).$$

This allows to conclude using \cref{eq:lowerregreteq}.
\medskip

\begin{proof}[Proof of Lemma~\ref{lemma:concentrationgeom}]
Note that the probability that the sum of $n$ \iid geometric variables of parameter $p$ are smaller than $kn/p$ is exactly the probability that the sum of $kn/p$ \iid Bernoulli variables are larger than $n$. We can then use the Chernoff bound on these $kn/p$ Bernoulli variables. The same reasoning also leads to the second inequality.
\end{proof}

\subsection{Proof of Theorem~\ref{thm:imperfectloss}}
\label{app:proofs-dynamic-imperfect}

This proof uses the block partitioning given by \cref{eq:defntnglobal} and the same transition matrix and reward as in \cref{app:proofs_dynamic_binary}. It relies on intermediate results given by Lemma~\ref{lemma:imperfect}. Its proof can be found below.

\begin{lemma}\label{lemma:imperfect}
For any $t\in [t_k+1, t_{k+1}]$,
\begin{enumerate}
\item $\pi_{t}^\imp(Q_{t}) \leq c^{t-t_k}\pi_{t_k}^\imp(Q_{t})$;
\item $\mathbb{E}\left[\ln(\pi_t^\imp(q)) \mid (t_n)_n , \pi_{t_k}^\imp, Q_t\neq q\right] \leq -(t-t_k)\frac{\delta^4}{2\gamma^2+4\delta^2} -\ln(\pi_{t_k}^\imp(Q_{t}))$;
\item $\mathbb{P}\Big[\ln(\pi_{t}^\imp(q)) - \mathbb{E}\left[\ln(\pi_{t}^\imp(q)) \mid (t_n)_n , \pi_{t_k}^\imp \right]  \geq \lambda \gamma d \sqrt{t-t_k} | (t_n)_n , \pi_{t_k}^\imp\Big]\leq \exp\left(-2\lambda^2\right)$;
\end{enumerate}
where $c = \max_{\pi, z, q,q'} \frac{G(z, \pi, q)}{G(z, \pi, q)}$.
\end{lemma}

Consider two successive blocks $k$ and $k+1$, where the quality is $q$ on the block $k$ and $q'$ on the block $k+1$. Similarly to the proof of Theorem~\ref{thm:dynamic-binary}, define $\tau_{k+1} = \min\Big(\left\{ t \in [t_{k+1}+1, t_{k+2}] \mid \pi_{t}^\imp(q') \geq 1/2 \right\} \cup \{ t_{k+2}\}\Big)$. We define $\tau_k$ similarly.

The first point of Lemma~\ref{lemma:imperfect} implies that $\tau_{k+1}-t_{k+1} \geq \min\Big(t_{k+2}-t_{k+1}, \frac{-\ln(2)-\ln(\pi^\imp_{t_{k+1}}(q'))}{\ln(c)}\Big)$.

\medskip

Moreover, thanks to the second and third points of Lemma~\ref{lemma:imperfect}, with probability at least $1-e^{-2\lambda^2}$ for some $\lambda >0$, \begin{equation}\label{eq:impln}
-\ln(\pi^\imp_{t_{k+1}}(q')) \geq (t_{k+1}-t_k)\frac{\delta^4}{2\gamma^2+4\delta^2} + \ln(\pi^\imp_{t_{k}}(q)) - \lambda \gamma d \sqrt{t_{k+1}-t_k}.
\end{equation}

Either $\pi^\imp_{t_{k}}(q) \geq  \exp\left(-(t_{k+1}-t_k)\frac{\delta^4}{4\gamma^2+8\delta^2} \right)$, in which case the two first terms in \cref{eq:impln} are larger than $(t_{k+1}-t_k)\frac{\delta^4}{4\gamma^2+8\delta^2}$.

Otherwise, $\pi^\imp_{t_{k}}(q) \leq  \exp\left(-(t_{k+1}-t_k)\frac{\delta^4}{4\gamma^2+8\delta^2} \right)$. Using the first point of \cref{lemma:imperfect}, this yields that for the $\frac{-\ln(2)}{\ln(c)}+(t_{k+1}-t_k)\frac{\delta^4}{(4\gamma^2+8\delta^2)\ln(c)}$ first steps of the block $k$, $\pi^\imp_{t}(q) \leq \frac{1}{2}$, \ie $\tau_k-t_k \geq \frac{-\ln(2)}{\ln(c)}+(t_{k+1}-t_k)\frac{\delta^4}{(4\gamma^2+8\delta^2)\ln(c)}$.

\medskip

So we can actually bound the error in expectation:

\begin{equation}\label{eq:boundtauimperfect}
\begin{aligned}
\mathbb{E}[\tau_{k+1} - t_{k+1} + \tau_k-t_k | (t_n)_n] \geq &  (1-e^{-2\lambda^2}) \min\left(t_{k+2}-t_{k+1}, \frac{\delta^4 \ (t_{k+1}-t_k)}{(4\gamma^2+8\delta^2)\max(1,\ln(c))}  \right) \\ & - \frac{\ln(2)+\lambda \gamma d\sqrt{t_{k+1}-t_k}}{\ln(c)}.
\end{aligned}
\end{equation}

\paragraph{Case $ \eta \geq 32$.}%
Consider $\beta \in [\frac{1}{4}, \frac{1}{2}]$ such that $\beta  \eta T \in 2\mathbb{N}^*$. As $t_{k+1}-t_k$ are dominated by geometric variables of parameter $ \eta$, we can show similarly to Equations~\eqref{eq:geombound1} and \eqref{eq:geombound2} in the proof of Theorem~\ref{thm:dynamic-binary} that
\begin{enumerate}
\item $\mathbb{P}\left[ \sum_{k=1}^{\beta \eta T} t_{k+1}-t_k > T \right] \leq e^{-4}$;
\item  $\mathbb{P}\left[ \sum_{k=1}^{\beta \eta T/2} \one_{t_{2k+1}-t_{2k} \geq \frac{2}{ \eta}} \one_{t_{2k+2}-t_{2k+1} \geq \frac{2}{\eta}} \leq \frac{\beta \eta T}{4} (1- \eta/2)^{\frac{4}{ \eta}} \right] \leq\exp\left(-\frac{7e^{-3}}{8}\right)$.
\end{enumerate}

Similarly to the proof of Theorem~\ref{thm:dynamic-binary}, the sum of these two probabilities is below $1$, so that none of these two events can happen with probability $\Omega(1)$. When it is the case, the first point yields that the $\beta \eta T$ first blocks totally count in the estimation error before $T$. The second point implies, thanks to Equation~\eqref{eq:boundtauimperfect}, that the estimation loss is $\Omega(T)$ in this case.

\paragraph{Case $ \eta T \leq 32$.}
Since $ \eta T = \Omega(1)$, we can consider a constant $c_0 >0$ such that $ \eta T > c_0$. Similarly to the case $ \eta T \leq 32$ in the proof of Theorem~\ref{thm:dynamic-binary}, we can show that with a positive probability depending only on $c_0$, the two first blocks are completed before $T$ and $\min(t_1 - t_0, t_2-t_1) \geq T/4$. In that case, Equation~\eqref{eq:boundtauimperfect} yields that the estimation loss incurred during the two first blocks is $\Omega(T)$, which leads to a regret $\Omega(MdT)$.

\bigskip

\begin{proof}[Proof of Lemma~\ref{lemma:imperfect}]\hfill \\
1) This is a direct consequence of the posterior update given by Equation~\eqref{eq:updatestationarydiscretefull}.\vspace{1em}\\
2) Jensen's inequality gives that 
$$\mathbb{E}\left[\ln(\pi_{t}^\imp(q)) \mid (t_n)_n , \pi_{t_k}^\imp\right] \leq \ln\left(\mathbb{E}\left[\pi_{t}^\imp(q) \mid (t_n)_n , \pi_{t_k}^\imp\right]\right).$$
Theorem~\ref{thm:stationary_discrete} claims that $$\mathbb{E}\left[\pi_{t}^\imp(q) \mid (t_n)_n , \pi_{t_k}^\imp\right] \leq\exp\left(-(t-t_k)\frac{\delta^4}{2\gamma^2+4\delta^2}\right)\frac{1}{\pi_{t_k}^\imp(Q_{t})},$$leading to the second point.\vspace{1em}\\
3) Recall that $\ln(\pi_{t}^\imp(q)) = \ln(\pi_{t_k}^\imp(q)) + \sum_{s=t_k}^{t-1}\ln\left(\frac{G(Z_s, \pi_s, q)}{G(Z_s, \pi_s)}\right)$ and that $\ln\left(\frac{G(Z_s, \pi_s, q)}{G(Z_s, \pi_s)}\right) \in [Y_s, Y_s + \gamma d]$ for some variable $Y_s$. The third point is then a direct application of Azuma-Hoeffding's inequality as used in the proof of Theorem~\ref{thm:stationary_discrete}.
\end{proof}

\section{Continuous quality} \label{app:continuous}

We consider in this section the continuous case where $\mathcal{Q}$ is some continuous set and show that, in the dynamic model described by \cref{eq:dynamic}, the regret is upper bounded by $\bigO{M\eta^{1/4}T}$ and lower bounded by $\Omega(M\eta^{1/2}T)$ when the reward function is $M$-Lipschitz. Closing the gap between these two bounds is left open for future work.

\subsection{Continuous model}

In the whole section, the quality space $\mathcal{Q}$ is a convex and compact subset of $\mathbb{R}^d$. \cref{ass:buyers} is specific to the discrete model and we use an equivalent assumption in the continuous case.

\begin{assumption}[Purchase guarantee, continuous case]\label{ass:buyerscontinuous}
The function $r$ is non-decreasing in each feature $q^{(i)}$ and there is some $\munderbar{q} \in \mathbb{R}^d$ such that $\forall i\in[d], q\in\mathcal{Q}, \munderbar{q}^{(i)} \leq q^{(i)}$ and $\mathbb{P}_{\theta_t}\big(r(\munderbar{q}, \theta_t) > 0 \big) > 0$.
\end{assumption}

In the continuous case, an additional assumption is required to get fast convergence of the posterior. 

\begin{assumption}[Monotone feedback]\label{ass:identification2}
For any $i \in \lbrace 1, \ldots, d \rbrace$ and $\pi_t \in \mathcal{P}(\mathcal{Q})$, $G^{(i)}(z^{(i)}, \pi_t, \cdot)$ defined by Equation~\eqref{eq:defnGi} is continuously differentiable and strictly monotone in $q^{(i)}$ for some $z \in \mathcal{Z}$.
\end{assumption}

This assumption guarantees that for two different qualities, the distributions of observed feedbacks are different enough. Note that $G_i$ does not have to be strictly monotone in $q^{(i)}$ for all $z \in \mathcal{Z}$, but only for one of them. For instance in the sparse feedback model, the probability of observing $z^{(i)}=\ast$ indeed does not depend on the quality as it corresponds to the absence of review. 
Requiring the monotonicity only for some $z_i$ is thus much weaker than for all of them.
\subsection{Stationary environment}\label{sec:stationary_continuous}

Consider as a warmup in this section the static case $Q_t = Q_1$ for all $t\in\mathbb{N}$.
The arguments from Section~\ref{sec:stationary} cannot be adapted to this case for two reasons. 
First, the pointwise convergence was shown using the fact that the posterior was upper bounded by $1$, but a similar bound does not hold for density functions.
Second, even the pointwise convergence of the posterior does not give a good enough rate of convergence for the estimated quality. 
Instead, we first show the existence of a ``good'' non-Bayesian estimator. 
The Bayes estimator will also  have similar, if not better, performances  as it minimizes the Bayesian risk.

\medskip

We first show the existence of a good non-Bayesian estimator. 
Define $L_t(z)$ as the empirical probability of observing the feedback $z$, \ie $L_t(z) = \frac{1}{t}\sum_{s=1}^{t-1} \one_{Z_s=z}$. Also define for any posterior $\pi$ and quality $q$:
\begin{equation}\label{eq:defnpsi}
\psi(\pi, q)  \coloneqq \left( z \mapsto G(z, \pi, q)\right),
\end{equation}
where $G$ is defined by Equation~\eqref{eq:defnG}. 
The function $\psi(\pi, q)$ is simply the probability distribution of the feedback, given the posterior $\pi$ and the quality $q$.
\begin{lemma}\label{lemma:estimrate}
Under \cref{ass:buyerscontinuous,ass:identification2},
\begin{align*}
\mathbb{E}\left[\left\|\bar{\psi}_{t+1}^{\dagger}\big(L_{t+1}
\big) - Q\right\|^2_2 \right] &= \bigO{1/t},
\shortintertext{where}
\bar{\psi}_{t+1}(\argdot) &\coloneqq \frac{1}{t}  {\textstyle\sum_{s=1}^{t} \psi(\pi_s, \argdot)} 
\shortintertext{and } 
\bar{\psi}_{t+1}^{\dagger}\big(L_{t+1}
\big) & \coloneqq \argmin_{Q \in \mathcal{Q}} \|L_{t+1}-\bar{\psi}_{t+1}(Q)\|^2_2 & \\  & = \argmin_{Q \in \mathcal{Q}} \sum_{z \in \mathcal{Z}}(L_{t+1}(z)-\bar{\psi}_{t+1}(Q)(z))^2.&
\end{align*}
\end{lemma}

The $\dagger$ operator is a generalized inverse operator, \ie $f^\dagger$ is the composition of $f^{-1}$ with the projection on the image of $f$. For a bijective function, it is then exactly its inverse. 

The $\argmin$ above is well defined by continuity of $\bar{\psi}_{t+1}$ and compactness of $\mathcal{Q}$. 
Assumption~\ref{ass:identification2} implies that $\bar{\psi}_{t+1}$ is injective. Thanks to this, the function $\bar{\psi}_{t+1}^\dagger$ is well defined.

Here $L_{t+1}$ is the empirical distribution of the feedback. The function $\bar{\psi}_{t+1}^{\dagger}$ then returns the quality that best fits this empirical distribution.

\begin{proof}
Note that $L_{t+1}(z)=\frac{1}{t}\sum_{s=1}^{t} \one_{Z_{s}=z}$, where $\mathbb{E}\left[\one_{Z_{s}=z} \mid \mathcal{H}_{s}, Q\right] = G(z, \pi_s, Q)$. As we consider the variance of a sum of martingales, we have

\begin{equation*}
\mathbb{E}\left[\left(L_{t+1}(z) - \bar{\psi}_{t+1}(Q)(z)  \right)^2 \ \Big| \ Q \right]  = \frac{1}{t^2}\sum_{s=1}^t \mathrm{Var}(\one_{Z_s=z} \mid Q, \pi_s).
\end{equation*}
From this, we deduce a convergence rate $1/t$:
\begin{align*}
\mathbb{E}\left[ \| L_{t+1} - \bar{\psi}_{t+1}(Q)\|_2^2  \,\Big|\, Q \right] 
 & \leq \frac{1}{t^2}\sum_{s=1}^t\sum_{z \in \mathcal{Z}} \mathrm{Var}(\one_{Z_s=z} \mid Q, \pi_s) \\
 & \leq \frac{1}{t^2}\sum_{s=1}^t\sum_{z \in \mathcal{Z}} \mathbb{P}(Z_s=z \mid Q, \pi_s) = \frac{1}{t}.\stepcounter{equation}\tag{\theequation}\label{eq:nonbayesian1}
\end{align*}

As $G^{(i)}$ is strictly monotone in $q^{(i)}$ and continuously differentiable on $\mathcal{Q}$ for some $z^{(i)}$, the absolute value of its derivative in $q^{(i)}$ is lower bounded by some positive constant. 
As a consequence, for some $\lambda > 0$,
\begin{equation}\label{eq:psiregular}
\forall q, q' \in \mathcal{Q}, \| q - q' \| \leq \lambda  \| \bar{\psi}_{t+1}(q) - \bar{\psi}_{t+1}(q') \|.
\end{equation}

For $\hat{Q} = \bar{\psi}_{t+1}^\dagger(L_{t+1}) = \argmin_{Q \in \mathcal{Q}} \|L_{t+1}-\bar{\psi}_{t+1}(Q)\|^2$, it follows
\begin{align*}
\mathbb{E}\left[ \| \hat{Q} - Q\|_2^2  \,\Big|\, Q \right] & \leq \lambda \mathbb{E}\left[ \| \bar{\psi}_{t+1}(\hat{Q}) - \bar{\psi}_{t+1}(Q)\|_2^2  \,\Big|\, Q \right]  \\
& \leq 2\lambda \mathbb{E}\left[ \| L_{t+1} - \bar{\psi}_{t+1}(\hat{Q})\|_2^2  \,\Big|\, Q \right] + 2\lambda \mathbb{E}\left[ \| L_{t+1} - \bar{\psi}_{t+1}(Q)\|_2^2  \,\Big|\, Q \right]  \\
& \leq 4 \lambda \mathbb{E}\left[ \| L_{t+1} - \bar{\psi}_{t+1}(Q)\|_2^2  \,\Big|\, Q \right] \leq  \frac{ \lambda}{t}.
\end{align*}
The third inequality is given by the definition of $\hat{Q}$ as a minimizer of the distance to $L_{t+1}$, and Lemma~\ref{lemma:estimrate} follows thanks to Equation~\eqref{eq:nonbayesian1}.
\end{proof}

In the sequel we use the notation $M_t = \mathbb{E}[Q \mid \mathcal{H}_t]$. 
Lemma~\ref{lemma:estimrate} gives a non-Bayesian estimator that converges to $Q$ at rate~$1/t$ in quadratic loss. Using arguments similar to \cite{besbes2018}, this implies that $M_t \xrightarrow{a.s.} Q$, thanks to a result from \cite{lecam2000}. Theorem~\ref{thm:continuousconvergence} yields a different result: $M_t$ converges to $Q$ at a rate $1/\sqrt{t}$ in average.
\begin{theorem} \label{thm:continuousconvergence}
Under \cref{ass:buyerscontinuous,ass:identification2}, then
$\mathbb{E}\Big[\left\|M_{t+1} - Q\right\|_1\Big] = \bigO{\sqrt{d/t}}$.
\end{theorem}
The hidden constant in the $\bigO{\cdot}$ above depends only the parameter $\lambda$ appearing in \cref{eq:psiregular}, which depends on the functions $G^{(i)}$. The above bound directly leads to a $\bigO{\sqrt{d/t}}$ regret when the reward is $M$-Lipschitz (for the $1$-norm).

Note that the rate $\bigO{\sqrt{d/t}}$ is the best rate possible even if the reviews report exactly $Q + \varepsilon_t$. Indeed, the best estimator in this case is the average quality $\frac{1}{t}\sum_{s<t} (Q+\varepsilon_s)$, which behaves as a Gaussian variable with a variance of order $1/t$ by the central limit theorem on each coordinate.
The error $\mathbb{E}\Big[\left\|M_{t+1} - Q\right\|\Big]$ is then of the same order as the square root of the trace of the covariance matrix of $M_{t+1}$, \ie $\sqrt{d/t}$.

\begin{proof}
A characterization of the Bayes estimator is that it minimizes the Bayesian mean square error among all $\mathcal{H}_t$-measurable functions. In particular, 
\begin{equation*}
\mathbb{E}\Big[\left\|M_{t+1} - Q\right\|^2_2 \Big] \leq \mathbb{E}\left[\left\|\bar{\psi}_{t+1}^{\dagger}\big(L_{t+1}
\big) - Q\right\|^2_ 2 \right].
\end{equation*}
Thanks to Lemma~\ref{lemma:estimrate}, this term is  $\bigO{1/t}$ and Theorem~\ref{thm:continuousconvergence} then follows by comparison of the $1$ and $2$-norms. 
\end{proof}

\medskip

\subsection{Dynamical environment}\label{sec:dynamic_continuous}

We now consider the dynamical setting given by Equation~\eqref{eq:dynamic}. The Markov chain is here continuous, but the quality still has a probability to stay the same $1-\eta$ at each round. As in the stationary case, we first expose a satisfying non-Bayesian estimator, implying similar bounds on the posterior distribution.
\medskip

In the stationary case, our non-Bayesian estimator comes from the empirical distribution of the feedback. As highlighted by Equation~\eqref{eq:posteriordynamic}, with a dynamical quality, recent reviews have a larger weight in the posterior. This leads to the following adapted discounted estimator for $\eta_1 \in (0,1)$:
\begin{equation}\label{eq:defnalpha}
L^{\eta_1}_t(z) \coloneqq \eta_1 {\textstyle\sum_{s=1}^{t-1}} (1-\eta_1)^{t-s-1} \one_{Z_{s}=z}.
\end{equation}

Lemma~\ref{lemma:estimratedynamic} below bounds the mean error for the estimator~$L^{\eta_1}_t$.
\begin{lemma}\label{lemma:estimratedynamic}
Under \cref{ass:buyerscontinuous,ass:identification2}, for $\eta_1 = \sqrt{\eta}$, 
\begin{equation*}
\sum_{t=1}^T \sqrt{\mathbb{E}\left[\left\|L^{\eta_1}_t(z)-\bar{\psi}_{t, \eta_1}(Q_{t})\right\|^2_2 \right]} = \bigO{\eta^{1/4} T},
\end{equation*}
where 
\begin{equation*}
\bar{\psi}_{t, \eta_1}(Q)(z) \coloneqq \eta_1 \sum_{s=1}^{t-1}(1-\eta_1)^{t-s-1} G(z, \pi_s, Q).
\end{equation*}
\end{lemma}

\begin{proof}
First fix the qualities $(Q_s)_s$ and blocks $(t_n)_n$ defined as in Equation~\eqref{eq:defntnglobal}. Note that $G(z, \pi_t, Q_{t})$ is exactly the expectation of $\one_{Z_{t}=z}$ given $\mathcal{H}_{t}$ and $Q_t$. Similarly to the stationary case, we have
\begin{align*}
\mathbb{E}\left[ \left(L^{\eta_1}_t(z) - \eta_1 \sum_{s=1}^{t-1}(1-\eta_1)^{t-s-1} G(z, \pi_s, Q_{s}) \right)^2 \ \bigg\vert \ (Q_s)_s \right] 
& =  \eta_1^2 \sum_{s=1}^{t-1}(1-\eta_1)^{2(t-s-1)} \mathrm{Var}(\one_{Z_s=z} \mid Q_s).
\end{align*}%

When summing over all $z \in \mathcal{Z}$, we get the following inequality
\begin{equation}\label{eq:concentrationdynamic1}
\mathbb{E}\left[ \Big\|L^{\eta_1}_t - \eta_1 \sum_{s=1}^{t-1}(1-\eta_1)^{t-s-1} G(\cdot, \pi_s, Q_{s}) \Big\|^2  \bigg| \ (Q_s)_s \right] \leq  \frac{\eta_1^2}{1-(1-\eta_ 1)^2} \leq \eta_1.
\end{equation}


For $t \in [t_i+1, t_{i+1}]$, we can relate the expected value of $L^{\eta_1}_t(z)$ to $\bar{\psi}_{t, \eta_1}(Q_t)(z)$:
\begin{align*}
\Big(\eta_1 \sum_{s=1}^{t-1}(1-\eta_1)^{t-s-1} G(z, \pi_s, Q_{s}) - \bar{\psi}_{t, \eta_1}(Q_t)(z)\Big)^2 & = \eta_1^2 \Big(\sum_{s=1}^{t-1}(1-\eta_1)^{t-s-1} \left(G(z, \pi_s, Q_{s})-G(z, \pi_s, Q_{t})\right)\Big)^2 \\
& = \eta_1^2 \Big(\sum_{s=1}^{t_i}(1-\eta_1)^{t-s-1} \left(G(z, \pi_s, Q_{s})-G(z, \pi_s, Q_{t})\right)\Big)^2 \\
& \leq (1-\eta_1)^{2t-2(t_i+1)}.\stepcounter{equation}\tag{\theequation}\label{eq:concentrationdynamic2}
\end{align*}%
The second equality holds because $Q_s = Q_t$ for $s > t_i$ by definition of the blocks. 
In the last inequality, we used the fact that $G$ has values in $[0,1]$, besides comparing the partial sum with $(1-\eta_1)^{t-(t_i+1)}/\eta_ 1$.
This finally gives, for $h(t) \coloneqq \max \left( \{ t' < t \mid Q_{t'} \neq Q_{t} \} \cup \{ 0 \} \right)$,
\begin{equation}\label{eq:concentrationdynamic4}
\mathbb{E}\left[\Big(\eta_1 \sum_{s=1}^{t-1}(1-\eta_1)^{t-s-1} G(z, \pi_s, Q_{s}) - \bar{\psi}_{t, \eta_1}(Q_t)(z)\Big)^2 \ \Big\vert \  (Q_s)_s \right] \leq (1-\eta_1)^{2(t-h(t)-1)}.
\end{equation}

When reversing the time, note that $t-h(t)-1$ is the minimum between a geometric variable of parameter $\eta$ and~${t-1}$. It follows
\begin{align*}
\mathbb{E}\left[\Big(\eta_1 \sum_{s=1}^{t-1}(1-\eta_1)^{t-s-1} G(z, \pi_s, Q_{s}) - \bar{\psi}_{t, \eta_1}(Q_t)(z)\Big)^2\right] & \leq \eta\sum_{s=0}^{\infty}(1-\eta)^s (1-\eta_1)^{2s} \\
& \leq \frac{\eta}{1-(1-\eta)(1-\eta_1)^2} \\
& \leq \frac{\eta}{\eta + \eta_1 - \eta_1 \eta}.\stepcounter{equation}\tag{\theequation}\label{eq:concentrationdynamic5}
\end{align*}

Noting that $2x^2 + 2y^2 \geq (x+y)^2$, we can now use Equations~\eqref{eq:concentrationdynamic1} and \eqref{eq:concentrationdynamic5} to bound the total error on a round:
\begin{align*}
\mathbb{E}\left[\left(L_t^{\eta_1}(z)-\bar{\psi}_{t, \eta_1}(Q_t)(z)\right)^2  \right] \leq 2\eta_1  + \frac{2\eta}{\eta + \eta_1 - \eta_1 \eta}.
\end{align*}%
The error on a single round is of order $\bigO{\eta_1 + \frac{\eta}{\eta + \eta_1} }$ in average; and for $\eta_1 = \sqrt{\eta}$, it is then $\bigO{\sqrt{\eta}}$ in average. Summing the square root of this term over all rounds finally yields Lemma~\ref{lemma:estimratedynamic}. 
\end{proof}

\begin{theorem}
\label{thm:bounddynamiccontinuous} 
If the reward is $M$-Lipschitz (for the $1$-norm), the regret of Bayesian consumers in the dynamical continuous case is bounded as $\regret = \bigO{M \sqrt{d}\eta^{1/4}T}$ under \cref{ass:buyerscontinuous,ass:identification2}.
\end{theorem}

As in the stationary case, the hidden constant in the $\bigO{\cdot}$ above depends only the parameter $\lambda$ appearing in \cref{eq:psiregular}.

\begin{proof}[Proof of Theorem~\ref{thm:bounddynamiccontinuous} ]
Similarly to the stationary setting, the error of the Bayesian estimator can be bounded by the error of the non-Bayesian one since the former is the minimizer of the quadratic loss among all $\mathcal{H}_t$-measurable functions:
\begin{align*}
\mathbb{E}\Big[\big\|\mathbb{E}\left[ Q_t \mid \mathcal{H}_t \right] - Q_t\big\|^2 \Big] & \leq \mathbb{E}\left[\left\| \bar{\psi}_{t, \eta_1}^{\dagger}(L^{\eta_1}_t)-Q_t\right\|^2 \right].
\end{align*}%
Thanks to Assumption~\ref{ass:identification2}, $\bar{\psi}_{t, \eta_1}$ verifies Equation~\eqref{eq:psiregular} for some constant $\lambda>0$ independent from $\eta_1$ and $T$ for any $t\geq \frac{1}{\eta_1}$. 
As we consider $\eta T = \Omega(1)$, the $\frac{1}{\eta_1}$ first terms in the loss are negligible compared $\eta^{1/4}T$.
The convergence rate is thus preserved when composing with $\bar{\psi}_{t, \eta_1}^{\dagger}$. Theorem~\ref{thm:bounddynamiccontinuous} then follows using Lemma~\ref{lemma:estimratedynamic} and Jensen's inequality as in the proof of Theorem~\ref{thm:continuousconvergence}.%
\end{proof}

\medskip

In contrast to the discrete case, determining a tight bound in the continuous case remains open for the dynamical setting. Note that the total error is of order at least $M\sqrt{d\eta}T$. Indeed, in the stationary case, no estimator converges faster than a rate $\sqrt{d/t}$. As the length of a block is around $1/\eta$, the loss per block is thus $\Omega\left(\sqrt{d/\eta}\right)$. 
Thanks to this, a tight bound should be between $M\sqrt{d\eta} T$ and $M\sqrt{d}\eta^{1/4} T$.

\medskip

A reason for such a discrepancy between the discrete and continuous case might be that the analysis is not tight enough. Especially, the considered non-Bayesian estimators might have a much larger regret than the Bayesian estimator.

On the other hand, reversing the posterior belief after a change of quality already takes a considerable amount of time in the discrete case and causes a loss $\ln(1/\eta)$ against $1$ in the stationary case. Here as well, reversing this belief might take a larger time, causing a loss $\eta^{-\frac{3}{4}}$ per block, against $\eta^{-\frac{1}{2}}$ in the stationary case. 
Showing a tighter lower bound is yet much harder than for the discrete case, as working directly on the Bayesian estimator is more intricate.

Lemma~\ref{lemma:estimratedynamic} uses the non-Bayesian estimator $L^{\eta_1}_t$ with the parameter $\eta_1$.
Quite surprisingly, $\sqrt{\eta}$ seems to be the best choice for the parameter $\eta_1$, despite $\eta$ being the natural choice. Figure~\ref{fullfig} below confirms this point empirically on a toy example.
The code used for this experiment can be found in the supplementary material.
The experiment considers the classical unidimensional setting with:
\begin{gather*}
r(Q, \theta) = Q + \theta, \\
f(Q, \theta, \varepsilon) = \mathrm{sign}(Q + \theta + \varepsilon).
\end{gather*}
Here, $\mathcal{Q}=[0,1]$, $\eta = 10^{-4}$ and $\theta$ and $\varepsilon$ both have Gaussian distributions. 
The Markov Chain is here given as follows
\begin{equation*}
\begin{cases}
Q_{t+1} = Q_t \text{ with probability } 1-\eta \\
Q_{t+1} = X_{t+1} \text{ otherwise},
\end{cases}
\end{equation*}
where $(X_t)$ is an i.i.d. sequence of random variables drawn from the uniform distribution on $[0,1]$.

Computing the exact posterior $M_t = \mathbb{E}\left[ Q_t \mid \mathcal{H}_t \right]$ is intractable, so we remedy this point by assuming $M_t = 1$ all the time. This simplification does not affect the experiments run here as $\bar{\psi}_{t, \eta_1}$ uses $M_t$ only to determine the population of potential buyers.

%
%
\begin{figure}[htbp]
\begin{adjustwidth}{-50pt}{-50pt}\centering
\caption{\label{fullfig}Behavior of $L^{\eta_1}$ for different $\eta_1$.}
\vspace{0.5em}
\subfigure[Estimation error of $L^{\eta_1}$. The error is $\sum_{t=1}^T  \sqrt{\mathbb{E}\left[\left\|L^{\eta_1}_t-\bar{\psi}_{t, \eta_1}(Q_t)\right\|_2^2 \right]}$ for ${T=10^5}$, where the expectation is estimated by averaging over $2000$ instances.]{\label{table:eta1}\centering \hspace{7em}
\def\arraystretch{1.5}
 \begin{tabular}{|x{55pt} || x{30pt}|x{30pt}|x{30pt}|x{30pt}|} 
 \hline
 \textbf{Value of $\pmb{\eta_1}$} & $\eta^{1/3}$ &  $\eta^{1/2}$ & $\eta^{2/3}$  & $\eta$ \\ 
 \hline
 \textbf{Error} & 10166 & 4780 & 3060 & 6462 \\ \hline
\end{tabular}
\hspace{7em}}\vspace{0.2em}
\subfigure[Tracking of $\bar{\psi}_{t, \eta_1}(Q_t)(1)$ by $L^{\eta_1}_t(1)$ over a single instance.]{\vspace{-0.5em}
\centering
\resizebox{\linewidth}{!}{\input{figure_eta1.pgf}}
}
\end{adjustwidth}
\end{figure}

A larger $\eta_1$ allows to forget faster past reviews and thus gives a better adaptation after a quality change. However, a larger $\eta_1$ also yields a less accurate estimator in stationary phases. 

The choice $\eta^{2/3}$ seems to be the best trade-off in Figure~\ref{fullfig}. The optimal choice of $\eta_1$ does not only depend on $\eta$ but also on the distributions of $\theta$ and $\varepsilon$. In the considered experiments, $\eta$ is thus not small enough to ignore these other dependencies. Figure~\ref{fullfig} yet illustrates the trade-off between small variance and fast adaptivity when tuning $\eta_1$.
\end{document}